\documentclass[12pt, dvips, twoside]{amsart}
\usepackage{amssymb,latexsym,bbm,graphicx,epsfig,epic,eepic,oldgerm}
\usepackage{psfrag}

\theoremstyle{plain}
\newtheorem{theorem}{Theorem}[section]
\newtheorem{lemma}[theorem]{Lemma}
\newtheorem{proposition}[theorem]{Proposition}
\newtheorem{corollary}[theorem]{Corollary}
\newtheorem{remark}[theorem]{Remark}
\newtheorem{remarks}[theorem]{Remarks}
\newtheorem{claim}[theorem]{Claim}

\newtheorem{definition}[theorem]{Definition}

\newtheorem{question}[theorem]{Question}

\newcommand{\proofend}{\hspace*{\fill} $\Box$\\}
\newcommand{\diam}{\hspace*{\fill} $\Diamond$\\}

\def\s{\smallskip}
\def\m{\medskip}
\def\eps{\varepsilon}
\def\Ker{\operatorname{ker}}
\def\Im{\operatorname {im}}

\def\Diff{\operatorname{Diff}}
\def\Diffcc{\operatorname{Diff^c}}
\def\Diffc0{\operatorname{Diff^c_0}}
\def\Symp{\operatorname{Symp}}
\def\Sympcc{\operatorname{Symp^c}}
\def\Sympc0{\operatorname{Symp^c_0}}

\def\Spin{\operatorname{Spin}}
\def\Int{\operatorname{Int}}
\def\Ham{\operatorname{Ham}}
\def\Hamc{\operatorname{Ham^c}}

\def\rank{\operatorname{rank}}

\def\length{\operatorname{length}}
\def\supp{\operatorname{supp}}
\def\idd{\operatorname{id}}
\def\ind{\operatorname{ind}}

\def\var{\operatorname{var}}
\def\can{\operatorname{can}}
\def\Flux{\operatorname{Flux}}

\def\reg{\operatorname{reg}}

\def\ga{\alpha}
\def\gb{\beta}
\def\gg{\gamma}
\def\gd{\delta}
\def\eps{\epsilon}

\def\gf{\varphi}

\def\gl{\lambda}
\def\go{\omega}
\def\gs{\sigma}
\def\gt{\vartheta}

\def\ca{{\mathcal A}}
\def\cb{{\mathcal B}}
\def\cc{{\mathcal C}}

\def\ce{{\mathcal E}}
\def\cf{{\mathcal F}}
\def\cg{{\mathcal G}}

\def\cj{{\mathcal J}}
\def\cl{{\mathcal L}}
\def\cm{{\mathcal M}}

\def\cp{{\mathcal P}}
\def\cR{{\mathcal R}}

\def\bJ{{\mathbf J}}
\def\bg{{\mathbf g}}

\def\CC{\mathbbm{C}}

\def\NN{\mathbbm{N}}

\def\RR{\mathbbm{R}}

\def\ZZ{\mathbbm{Z}}

\def\RP{\operatorname{\mathbbm{R}P}}  
\def\CP{\operatorname{\mathbbm{C}P}}
\def\HP{\operatorname{\mathbbm{H}P}}
\def\Ca{\operatorname{\mathbbm{C}\mathbbm{a}P}^2}

\def\Diff{\operatorname{Diff}}

\def\pp{\partial}

\def\ra{\rightarrow}
\def\ha{\hookrightarrow}

\def\ni{\noindent}
\def\b{\bigskip}
\def\m{\medskip}

\def\id{\mbox{id}}

\def\proof{\noindent {\it Proof. \;}}

\begin{document}

\title{Slow entropy and symplectomorphisms of cotangent bundles}    
$   $ \\
$   $ \\

\author{Urs Frauenfelder}
\address{(U.\ Frauenfelder) Department of Mathematics, Hokkaido University,
Sapporo 060-0810, Japan}
\email{urs@math.sci.hokudai.ac.jp}
\author{Felix Schlenk}
\address{(F.\ Schlenk) Mathematisches Institut,
Universit\"at Leipzig, 04109 Leipzig, Germany}
\email{schlenk@math.uni-leipzig.de}

\date{\today}

\begin{abstract}
We consider an entropy-type invariant which measures the polynomial
volume growth of submanifolds under the iterates of a map, and
we establish sharp
uniform lower bounds of this invariant for the following classes of
symplectomorphisms of cotangent bundles over a compact base:

\m
\begin{itemize}
\item[$\bullet$]
non-identical compactly supported symplectomorphisms which are
symplectically isotopic to the identity,
\item[$\bullet$]
symplectomorphisms generated by classical Hamiltonian functions,
\item[$\bullet$]
Dehn twist like symplectomorphisms  over compact rank one symmetric spaces.
\end{itemize}
\end{abstract}

\maketitle

\markboth{{\rm Slow entropy and symplectomorphisms of cotangent bundles}}{{}} 

\tableofcontents

\section{Introduction and main results}  

\ni
The topological entropy $h_{\text{top}} (\gf)$ of a compactly supported
smooth diffeomorphism $\gf$ of a smooth manifold $M$ is a basic numerical
invariant measuring the orbit structure complexity of $\gf$.
There are various ways of defining $h_{\text{top}} (\gf)$, see \cite{KH}.
A geometric way was found by Yomdin and Newhouse in their seminal works
\cite{Y} and \cite{N}:
Fix a Riemannian metric $g$ on $M$. 
For $i \in \left\{ 1, \dots, \dim M \right\}$ denote by $\Sigma_i$ 
the set of smooth embeddings $\gs$ of the cube $Q^i = [0,1]^i$
into $M$, and by $\mu_g(\gs)$ the volume of $\gs \left( Q^i \right) \subset M$
computed with respect to the measure on $\gs \left( Q^i \right)$
induced by $g$. 
The {\it $i$-dimensional volume growth}\, $v_i(\gf)$ of $\gf$ is defined as 
\[
v_i (\gf) \,=\, 
 \sup_{\gs \in \Sigma_i} \liminf_{n \ra \infty} 
  \frac{ \log \mu_g \left( \gf^n ( \gs ) \right)}{n} .
\]
Since $\gf$ is compactly supported,  
$v_i (\gf)$ does not depend on the choice of $g$ and is finite,
and $v_{\dim M} (\gf) = 0$.
The main result of \cite{Y,N} is
\begin{equation}  \label{e:ny}
\max_i v_i (\gf) \,=\,  h_{\text{top}} (\gf) .
\end{equation}

The purpose of this work is to study the volume growth of
symplectomorphisms of cotangent bundles $T^*B$ over a closed base $B$
endowed with their canonical symplectic structure $\go = d \gl$.
Cotangent bundles are the phase spaces of classical mechanics, and
classical Hamiltonian systems on such manifolds describe 
systems without friction.
The orbit structure of their time-1-maps is therefore often intricate,
and so one can expect that there exist non-trivial lower bounds of
entropy-type invariants for non-identical Hamiltonian diffeomorphisms on
cotangent bundles.

The topological entropy itself, which by \eqref{e:ny} measures 
the exponential volume growth, 
vanishes for many non-identical Hamiltonian diffeomorphisms.
Following \cite{KT}, we therefore look at the polynomial range of growth
and define for each compactly supported symplectomorphisms of $T^*B$ and each 
$i \in \left\{ 1, \dots, 2d = \dim T^*B \right\}$ 
the {\it $i$-dimensional slow volume growth}\, 
$s_i(\gf) \in [0, \infty]$ by
\begin{equation} \label{def:si}
s_i (\gf) \,=\, 
 \sup_{\gs \in \Sigma_i} \liminf_{n \ra \infty} 
    \frac{ \log \mu_g \left( \gf^n ( \gs ) \right)}{\log n}.
\end{equation}
Again, $s_i (\gf)$ does not depend on the choice of $g$, and
$s_{2d} (\gf) =0$.

We shall establish uniform lower bounds of $s_1$ for certain classes of 
symplectomorphisms which are symplectically isotopic to the identity 
and a uniform lower bound of $s_d$ for certain symplectomorphisms
some of which are smoothly but not symplectically isotopic to the identity.
For the moment, all Hamiltonian functions and all symplectomorphisms
are assumed to be $C^\infty$-smooth. 
Weaker smoothness assumptions are discussed in Section~\ref{smooth}.

We denote by $\Symp_0^c \left( T^*B \right)$ the identity component of
the group $\Symp^c \left( T^*B \right)$  of compactly supported
symplectomorphisms of $\left( T^*B, d \gl \right)$.
It contains the group of compactly supported Hamiltonian diffeomorphisms
generated by time-dependent compactly supported Hamiltonian functions $H
\colon [0,1] \times T^*B \ra \RR$.

\b
\ni
{\bf Theorem 1.}
{\it
For every non-identical symplectomorphism  
$\gf \in \Symp_0^c \left( T^*B \right)$ 
it holds true that $s_1(\gf) \ge 1$.
}

\b
\ni
Theorem 1 is sharp. Indeed, we shall show by means of an example
that

\b
\ni
{\bf Proposition 1.}
{\it
On every $2d$-dimensional symplectic manifold $(M, \go)$ there exists a
compactly supported Hamiltonian diffeomorphism $\gf$ such that $s_i(\gf)
= 1$ for all $i = 1, \dots, 2d-1$.
}

\b
\ni
Theorem~1 implies at once that the group $\Symp^c_0 \left( T^*B \right)$
has no torsion.
Stronger implications for the algebraic structure of the groups
$\Ham^c \left( T^*B \right)$ and $\Symp^c_0 \left( T^*B \right)$ 
are given in Section~\ref{s:t1}, where we work with arbitrary exact
symplectic manifolds convex at infinity.

\b
We next look at classical Hamiltonian systems.
We choose a Riemannian metric $g$ on $B$ and denote by $g^*$ the 
Riemannian metric induced on $T^*B$. 
We denote canonical coordinates on $T^*B$ by $(q,p)$. 
A classical Hamiltonian function $H \colon \RR \times T^*B \ra \RR$ 
is of the form
\[
H(t,q,p) \,=\, \tfrac{1}{2} \left| p - A(t,q) \right|^2 + V(t,q) 
\]
and periodic in the time variable $t$.
For the purpose of this paper we can assume without loss of generality
that the period is $1$. 
Writing $S^1 = \RR/ \ZZ$ we then have $H \colon S^1 \times T^*B \ra \RR$. 
It is well-known that such Hamiltonian functions generate a flow. We
denote its time-$1$-map by $\gf_H$, and we use the Riemannian metric
$g^*$ to define $s_1 \left( \gf_H \right)$ by \eqref{def:si}.
Since $\gf_H$ is not compactly supported, $s_1 \left( \gf_H \right)$
depends on the choice of $g$. 
For $r>0$ we abbreviate
\[
T_r^*B \,=\, \left\{ (q,p) \in T^*B \mid \left| p \right| \le r \right\}.
\]

\ni
{\bf Theorem~2.}
{\it
Assume that $B$ is a closed manifold whose fundamental group is
finite or contains infinitely many conjugacy classes, and that
$H \colon S^1 \times T^*B \ra \RR$
is a classical Hamiltonian function.
Then the following assertions hold true.

\begin{itemize}
\item[(i)]
$s_1(\gf_H) \ge 1/2$;
\item[(ii)]
$s_1(\gf_H) \ge 1$ provided that there exists $r < \infty$ such that
$H$ is time-independent on $T^*B \setminus T_r^*B$.
\end{itemize}
}

\b
\ni
Notice that the fundamental group $\pi_1(B)$ of a closed
manifold $B$ is finitely presented. 
As was pointed out to us by Indira Chatterji, Rostislav Grigorchuk
and Guido Mislin,
no infinite finitely presented group with finitely many conjugacy
classes is known, and so Theorem~2
possibly holds for all closed manifolds.

\b
\ni
{\bf Examples.}
We give two classes of closed manifolds $B$ for which the
assumption on $\pi_1(B)$ in Theorem~2 is met.
Denote by $\cc (B)$ the set of conjugacy classes of
$\pi_1(B)$.

\begin{itemize}
\item[1.]
Assume that $\pi_1(B)$ is abelian. Then $\cc (B) = \pi_1(B)$.
Examples of closed manifolds with abelian fundamental group are
Lie groups and, more generally, $H$-spaces.
\item[2.]
Assume that the first Betti number $b_1(B)$ does not
vanish. Then $\cc (B)$ is infinite. Indeed, the Hurewicz map
factors as 
\[
\pi_1(B) \,\ra\, \cc (B) \,\ra\, 
\frac{\pi_1(B)}{[\pi_1(B),\pi_1(B)]} = H_1(B;\ZZ) .
\]
\end{itemize}
We in particular see that Theorem~2 applies to all closed
$2$-manifolds and their products.
\diam

We finally look at certain compactly supported symplectomorphisms which
are not Hamiltonian.
The spaces we shall consider are the cotangent bundles over compact rank
one symmetric spaces (CROSS'es, for short), and the diffeomorphisms are
Dehn twist like symplectomorphisms.
These maps were introduced to symplectic topology by Arnol'd
\cite{A} and Seidel \cite{S1,S2}. 
They play a prominent role in the study of the symplectic
mapping class group of various symplectic manifolds, 
\cite{KS,S0,S1,S2,S8},
and generalized Dehn twists along spheres can be used to detect
symplectically knotted Lagrangian spheres, \cite{S1,S2},
and (partly through their appearance in Seidel's long exact
sequence in symplectic Floer homology) are an important
ingredient in attempts to prove Kontsevich's homological mirror
symmetry conjecture, \cite{KS,S4,S5,S7,S9,ST}.

\begin{figure}[h] 
 \begin{center}
  \psfrag{Tg}{$T^* \gg$}
  \psfrag{t}{$\gt |_{T^* \gg}$}
  \psfrag{g}{$\gg$}
  \leavevmode\epsfbox{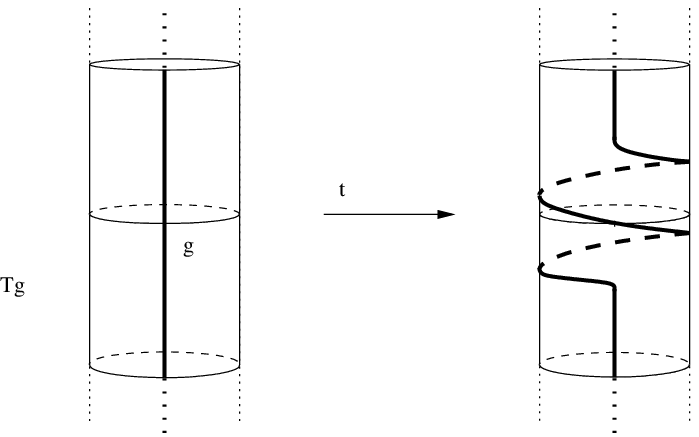}
 \end{center}
 \caption{The map $\gt |_{T^* \gg}$.}
 \label{figure1}
\end{figure}

Let $(B,g)$ be a CROSS, i.e., $B$ is a sphere $S^d$, a projective space
$\RP^d$ $\CP^n$, $\HP^n$, or the exceptional symmetric space 
$F_4 / \Spin_9$ diffeomorphic to the Cayley plane $\Ca$.
All geodesics on $(B,g)$ are embedded circles of equal length.
We define $\gt$ to be the compactly supported 
diffeomorphism of $T^*B$ whose restriction
to the cotangent bundle $T^* \gg \subset T^*B$ over any 
geodesic circle $\gg \subset B$ is the square of the ordinary Dehn twist along
$\gg$ depicted in Figure \ref{figure1}.

%
%

\ni
We call $\gt$ a {\it twist}. 
A more analytic description of twists is given in
Subsection~\ref{twists}.
It is known that twists are symplectic, and that the class of a twist
generates an infinite cyclic subgroup of the mapping class group
$\pi_0 \left( \Symp^c \left( T^*B \right)\right)$, see \cite{S2}.

A $d$-dimensional submanifold $L$ of $T^*B$
is called {\it Lagrangian}\, if $\omega$ vanishes on $TL \times TL$.
Lagrangian submanifolds play a fundamental role in symplectic geometry. 
For each $\gf \in \Symp^c \left( T^*B \right)$
we therefore also consider its {\it Lagrangian volume growth} 
\begin{equation*}  
l (\gf) \,=\, 
 \sup_{\gs \in \Lambda} \liminf_{n \ra \infty} 
    \frac{ \log \mu_g \left( \gf^n ( \gs ) \right)}{\log n}
\end{equation*}
where $\Lambda$ is the set smooth embeddings $\gs \colon Q^d \ha T^*B$
for which $\gs \left( Q^d \right)$ is a Lagrangian submanifold of $T^*B$.
Of course, $l (\gf) \le s_d (\gf)$.
As we shall see in Subsection~\ref{twists}, 
$s_i(\gt^m) = l(\gt^m) =1$ for every $i \in
\left\{1, \dots, 2d-1 \right\}$ and every $m \in \ZZ \setminus \{ 0\}$.

\b
\ni
{\bf Theorem 3}\:
{\it 
Let $B$ be a $d$-dimensional compact rank one symmetric space, and
let $\gt$ be the twist of $T^*B$ described above. 
Assume that $\gf \in \Symp^c \left( T^*B \right)$ 
is such that $[\gf] = \left[ \gt^m \right] 
\in \pi_0 \left( \Symp^c \left( T^*B \right) \right)$
for some $m \in \ZZ \setminus \{ 0 \}$. 
Then $s_d(\gf) \ge l(\gf) \ge 1$.
}

\b
Theorem~3 is of particular interest if $B$ is $S^{2n}$ or $\CP^n$, 
$n \ge 1$,
since in these cases it is known that $\gt$ can be deformed to the
identity through compactly supported {\it diffeomorphisms}, see
\cite{K, S2}.
Twists can be defined on the cotangent bundle of any Riemannian manifold
with periodic geodesic flow. 
In Section \ref{t3}, Theorem~3 is proved for all known such manifolds.

\begin{figure}[h] 
 \begin{center}
  \psfrag{Tg}{$T^* \gg$}
  \psfrag{t}{$\tau |_{T^* \gg}$}
  \psfrag{g}{$\gg$}
  \leavevmode\epsfbox{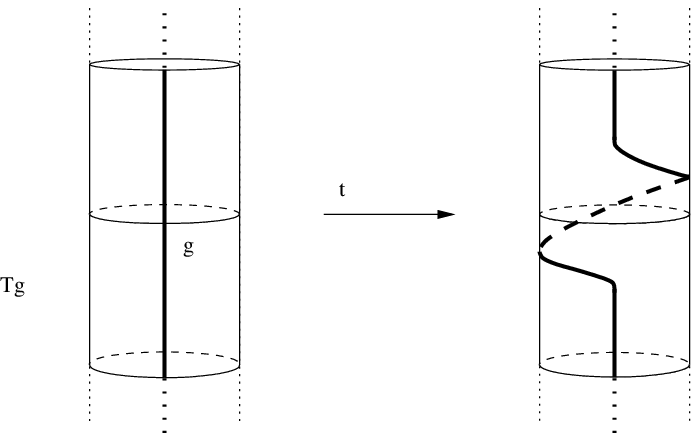}
 \end{center}
 \caption{The map $\tau |_{T^* \gg}$.}
 \label{figure2}
\end{figure}

In the case that $B$ is a sphere $S^d$, one can use the fact that all
geodesics emanating from a point meet again in the antipode to see
that the twist $\gt$ admits a square
root $\tau \in \Symp^c \left( T^*S^d \right)$. 
For $d = 1$, $\tau$ is the ordinary
Dehn twist along a circle, and for $d \ge 2$ it is the generalized
Dehn twist thoroughly studied in \cite{S0,S1,S2,S8}.
Given any great circle $\gg$ in $S^d$, the restriction of $\tau$ to
$T^* \gg \subset T^*S^d$ is the ordinary Dehn twist
along $\gg$ depicted in Figure~\ref{figure2}.

%
%

\b
\ni
{\bf Corollary 3}\:
{\it
Let $\tau$ be the (generalized) Dehn twist of $T^*S^d$ described above,
and assume that 
$\gf \in \Symp^c \left( T^*S^d \right)$ is 
such that $[\gf] = \left[ \tau^m \right]  
\in \pi_0 \left( \Sympcc \left( T^*S^d \right) \right)$
for some $m \in \ZZ \setminus \{ 0 \}$. 
Then $s_d(\gf) \ge l(\gf) \ge 1$. 
}

\b
\ni
Since $[\tau^2] = [\gt]$ has infinite order in 
$\pi_0 \left( \Sympcc \left( T^*S^d \right) \right)$, 
so has $[\tau]$.
In the case $d=2$ it is known that $[\tau]$ generates 
$\pi_0 \left( \Symp^c \left( T^*S^2 \right) \right)$, see \cite{S0,S8}.
Theorem~1 and Corollary~3 
thus give a nontrivial uniform lower bound of the slow volume growth
\[
s(\gf) \,=\, \max_i s_i(\gf)
\] 
for each $\gf \in \Sympcc (T^*S^2) \setminus \{ \id \}$.  

Following Shub \cite{Sh}, we consider a symplectomorphism $\gf \in
\Symp^c \left(T^*B \right)$ 
as a best diffeomorphism in its symplectic isotopy class if $\gf$
minimizes both
$s(\gf)$ and $l(\gf)$.
We shall show that 
$s_i(\tau^m) = l(\tau^m) =1$ and
$s_i(\gt^m) = l(\gt^m) =1$ for every $i \in
\left\{1, \dots, 2d-1 \right\}$ and every $m \in \ZZ \setminus \{ 0\}$.
In view of Theorem 3 and Corollary 3, the twists $\tau^m$ and
$\gt^m$ are then best diffeomorphisms in their symplectic isotopy classes.  


\m
Uniform lower bounds of an entropy type invariant 
were first obtained in a beautiful paper of Polterovich \cite{P1} 
for a class of symplectomorphisms of closed symplectic manifolds 
with vanishing second homotopy group.
E.g., it is shown that for any non-identical Hamiltonian diffeomorphism
$\gf$ of the standard $2d$-dimensional torus,
\begin{equation*}  
s_1 (\gf) \,\ge\,
 \left\{
  \begin{array}{rll} 
     1            &\text{ if} & d =1, \\
     \frac{1}{2}  &\text{ if} & d \ge 2. 
  \end{array}
 \right.
\end{equation*} 
The results in \cite{P1} are formulated in terms of the growth of the
uniform norm of the differential of $\gf$. In Section \ref{comp} we
shall reformulate our results in terms of this invariant.


\b
\ni
{\bf Acknowledgements.}
This paper owes very much to Leonid Polterovich and Paul Seidel.
We cordially thank both of them for their help.
We are grateful to Kai Cieliebak, Yuri Chekanov, 
Misha Entov, Janko Latschev, Anatole Katok, Dusa Mc\;\!Duff,
Dietmar Salamon, Matthias Schwarz and Edi Zehnder for useful discussions.
Much of this work was done during the stays of both authors at the
Symplectic Topology Program at Tel Aviv University in spring
2002
and at Hokkaido University at Sapporo in November 2003, 
and during the second authors stay at FIM of ETH Z\"urich and at
Leipzig University. 
We wish to thank these institutions for their kind hospitality, 
and the Swiss National Foundation for its generous support.

\section{Outline of the proofs}  \label{idea}

\ni 
As in the previous section we consider a smooth closed manifold $B$ and let
$T^*B$ be the cotangent bundle over $B$ endowed with its canonical 
symplectic form 
$\go \equiv \go_\text{can} = d \gl$ where 
$\gl = \sum p_i dq_i$.
We start with explaining our results for Hamiltonian diffeomorphisms.
The Hamiltonian diffeomorphisms addressed in Theorem~1 are end points of
smooth paths $\left\{ \gf^t \right\}$, $t \in [0,1]$, generated by
compactly supported Hamiltonian functions $H \colon [0,1] \times T^*B
\ra \RR$.
After reparametrizing the path $\left\{ \gf^t \right\}$ we can assume
that $H(t,x) = 0$ for $t$ near $0$ and $t$ near $1$, and so we can
assume that $H$ is defined on $S^1 \times T^*B$.
Such a Hamiltonian function or a classical Hamiltonian function generates
a flow $\gf_H^t$. 
We denote the set of its $1$-periodic orbits by $\cp_H$.
For $x \in \cp_H$ the symplectic action is defined as
\begin{equation}  \label{def:action}
\ca_H(x) \,=\, \int_0^1 x^*\gl - \int_0^1 H(t, x(t)) \, dt . 
\end{equation}
Assume that there exist $x,y \in \cp_H$ such that 
$\ca_H(x) < \ca_H(y)$. 
Following an idea of Polterovich \cite{P1},
we shall choose a smoothly embedded curve $\gs \colon [0,1] \ra T^*B$ 
connecting $x(0)$ with $y(0)$ and shall show that the quantity 
$\int_{\gf^n \gs} \gl$ grows linearly. 
From this we shall easily obtain  
\begin{itemize}
\item[(i)]
$s_1(\gf_H) \ge 1/2$;
\item[(ii)]
$s_1(\gf_H) \ge 1$ provided that there exists $r < \infty$ such that
\[
\left\{ \gf^n \gs \mid n \ge 1 \right\} \,\subset\, T_r^* B .
\]
\end{itemize}

\ni
We are then left with finding $x,y \in \cp_H$ such that
$\ca_H(x) < \ca_H(y)$. 
In the case of compactly supported Hamiltonian functions as in
Theorem~1, we shall do this by using a result from \cite{FS} 
relying on symplectic Floer homology.
In the case of the Hamiltonian functions considered in Theorem 2, we
shall use work of Benci \cite{Be} to show that the symplectic action
functional $\ca_H$ is not bounded from above on $\cp_H$.

\s 
Proving Theorem~1 for the whole group $\Symp_0^c \left( T^*B
\right)$ is now elementary: 
If $\dim B \ge 2$, every $\gf \in \Symp_0^c \left( T^*B
\right)$ is Hamiltonian, see the proof of Lemma~\ref{l:Ham} below;
and for a symplectomorphism $\gf \in \Symp_0^c \left( T^*S^1
\right)$ which is not Hamiltonian, the flux 
does not vanish, and this yields $s_1(\gf) \ge 1$ at once.

\s
The proof of Theorem~3 is different in nature.
To fix the ideas, we assume $B=S^d$, and that $\gf$ is isotopic to $\gt$
through symplectomorphisms supported in $T^*_1S^d$.
For $x \in S^d$ we denote by $D_x$ the $1$-disc in $T_x^*S^d$.
Consider first the case $d=1$. We fix $x$.
For a twist $\gt$ as in Figure~\ref{figure1} and $n \ge1$, the image
$\gt^n (D_x)$ wraps $2n$ times around the base $S^1$.
For topological reasons the same must hold for $\gf$, and so
\[
\mu_{g^*} \left( \gf^n (D_x) \right) \,\ge\, 2n \mu_g (S^1) .
\]
In particular, $s_1(\gf) \ge 1$.
For odd-dimensional spheres, Theorem~3 follows from a similar argument.
For even-dimensional spheres, however, Theorem~3 cannot hold for
topological reasons, since then $\gt$ is isotopic to the identity through
compactly supported diffeomorphisms.
In order to find a symplectic argument, we rephrase the above proof for
$S^1$ in symplectic terms: 
For every $y \neq x$ the Lagrangian submanifold $\gt^n(D_x)$ intersects
the Lagrangian submanifold $D_y$ in $2n$ points, and under symplectic
deformations of $\gt$ these $2n$ Lagrangian intersections persist. This
symplectic point of view generalizes to even-dimensional spheres:
For a twist $\gt$ on $S^d$ as in Figure~\ref{figure1}, $n \ge 1$
and $y \neq x$, the Lagrangian submanifolds $\gt^n (D_x)$ and $D_y$
intersect in exactly $2n$ points.
We shall prove that the Lagrangian Floer homology of $\gt^n (D_x)$ and $D_y$
has rank $2n$.
The isotopy invariance of Floer homology then implies that 
$\gf^n (D_x)$ and $D_y$ must intersect in at least $2n$ points.
Since this holds true for every $y \neq x$, we conclude that
\[
\mu_{g^*} \left( \gf^n (D_x) \right) \,\ge\, 2n \mu_g (S^d) .
\]
In particular, $s_d(\gf) \ge l(\gf) \ge 1$.

\b
For any compactly supported diffeomorphism $\gf$ of a manifold $M$
we denote by 
\[
\rho (\gf_*) \,=\, \lim_{n \ra \infty} \left\| \gf_*^n \right\|^{1/n}
\]
the spectral radius of the induced
automorphism $\gf_* \colon H_* (M;\RR) \ra H_*(M;\RR)$ of the total real
homology of $M$.
It only depends on the isotopy class of $\gf$.
On his search for the simplest diffeomorphism in each isotopy class, 
Shub \cite{Sh} formulated the entropy conjecture  
\begin{equation}  \label{est:ent}
h_\text{top}(\gf) \,\ge\, \log \rho (\gf_*) 
\end{equation}
for $C^1$-diffeomorphisms.
For $C^\infty$-diffeomorphisms the estimate
\eqref{est:ent} follows at once from Yomdin's estimate
$h_\text{top}(\gf) \ge \max_i v_i(\gf)$.
Besides certain Dehn twists 
all symplectomorphisms studied in this paper are isotopic to the
identity mapping, and so the estimate \eqref{est:ent} becomes vacuous.
Nevertheless, our results are of the same nature as the 
estimate~\eqref{est:ent}:
While the dynamical quantity $h_\text{top}(\gf)$ is replaced by the slow
volume growths $s_1(\gf)$ or $s_d(\gf)$ and $l(\gf)$, the homological
quantity $\log \rho (\gf_*)$ is replaced by
{\it Floer-homological}\, quantities.
In Theorems 1 and 2, this quantity is essentially the polynomial
growth rate 
of the difference of the
symplectic action of two closed orbits which represent generators of the
symplectic Floer homology of $\gf_H$,
and in Theorem 3 this quantity is the polynomial growth rate of the 
rank of a Lagrangian Floer homology associated with $\gf$.
In the case of odd-dimensional spheres, the lower bound $1$ in Theorem~3
can also be obtained by computing the homological growth of $\gt$, and
so in this case Theorem~3 is just a version of \eqref{est:ent}, 
see Subsection~\ref{var}.

\section{Slow entropy on $\Symp_0^c \left(M, d\gl \right)$}  \label{s:t1}
            
\ni
We consider a symplectic manifold $(M,\go)$ with or without boundary $\pp M$.
Denote by $\Symp^c_0 (M,\go)$ the identity component of the group 
$\Symp^c (M,\go)$ of symplectomorphisms of $(M, \go)$ whose support is
compact and contained in $M \setminus \pp M$,
and denote by $\Ham^c (M, \go)$ its subgroup consisting of
symplectomorphisms generated by 
Hamiltonian functions $H \colon S^1 \times M \ra \RR$ 
whose support is compact and contained in
$S^1 \times \left( M \setminus \pp M \right)$.
We shall only consider such Hamiltonians.
We define the slow volume growth $s_1(\gf)$ of $\gf \in
\Symp^c(M, \go)$ as in \eqref{def:si}.
A symplectic manifold $(M, \go)$ is called {\it exact}\, 
if there exists a $1$-form $\gl$ on $M$ such that $\go = d \gl$.
The boundary $\pp M$ of a $2d$-dimensional symplectic manifold 
$(M,\go)$ is said to be {\it convex}\,
if there exists a Liouville vector field $Y$ 
(i.e., $\cl_Y \go = d \iota_Y \go = \go$) 
which is defined near $\pp M$ and is everywhere transverse to $\pp M$, 
pointing outwards.
Equivalently, there exists a $1$-form $\ga$ on $\pp M$ such that $d \ga
= \go |_{\pp M}$ and such that $\ga \wedge (d \ga)^{d-1}$ is a volume
form inducing the boundary orientation of $\pp M \subset M$.
Following \cite{EG} we say that a symplectic manifold $(M, \go)$ is {\it
convex at infinity}\, if there exists an increasing sequence of 
compact submanifolds $M_i \subset M$ with smooth convex boundaries 
exhausting $M$, that is,
\[
M_1 \subset M_2 \subset \dots \subset M_i \subset \dots \subset M
\quad \text{ and } \quad
\bigcup_i M_i = M . 
\]
Cotangent bundles $T^*B$ over a closed base $B$ are exhausted by the
compact submanifolds $T_i^*B$ with convex boundary, and so 
Theorem~1 follows from 
\begin{theorem}  \label{t:4exact}
Assume that $(M, d\gl)$ is an exact symplectic manifold 
and that $\gf \in \Symp_0^c (M, d\gl) \setminus \left\{ \idd \right\}$.
If $\gf$ is Hamiltonian, assume in addition that $(M, d\gl)$ is
convex at infinity.
Then $s_1(\gf) \ge 1$.
\end{theorem}

\subsection{Proof of Theorem~1}  \label{ss:t1}
The main ingredient in the proof of Theorem~\ref{t:4exact} is 
the following result whose proof in \cite{FS} relies on
symplectic Floer homology.

\begin{proposition}  \label{p:4axa}
Assume that $(M, d\gl)$ is an exact symplectic manifold convex
at infinity. 
For any Hamiltonian function $H \colon S^1 \times M \ra \RR$ generating 
$\gf \in \Ham^c (M, d\gl) \setminus \left\{ \idd \right\}$ there exist
$x,y \in \cp_H$ such that $\ca_H(x) \neq \ca_H(y)$.
\end{proposition}

Here, $\cp_H$ denotes the set of $1$-periodic orbits of
$\gf_H^t$, and the symplectic action $\ca_H(x)$ of $x \in \cp_H$
is defined as in \eqref{def:action}.
In the remainder of the proof we closely follow the proof of 
Theorem~1.4.A in \cite{P1}.
Consider an exact symplectic manifold $(M, d \gl)$. If $\pp M$ is not
empty, we replace $M$ by its interior $M \setminus \pp M$, which we
denote again by $M$.
According to \cite[Chapter~10]{MS}, the map
\begin{equation}  \label{def:flux}
\Flux \colon \Symp_0^c (M) \ra H_c^1(M;\RR), \qquad \gf \mapsto 
[ \gf^* \gl - \gl] ,
\end{equation}
is a homomorphism which fits into the exact sequence
\begin{equation}  \label{s:exact}
0 \,\ra\, \Hamc (M) \,\ra\, \Sympc0 (M) \,\ra\, H_c^1 ( M;\RR) \,\ra\, 0.
\end{equation} 

\s
\ni
{\bf Case 1. $\gf$ is Hamiltonian.}
Assume that $\gf$ is generated by $H \colon S^1 \times M \ra
\RR$. According to Proposition~\ref{p:4axa} 
we find $x,y \in \cp_H$ such that 
\[
c \,:=\, \ca_H(y) - \ca_H(x) \,>\, 0 .
\] 
Choose a smoothly embedded curve $\gs \colon [0,1] \ra M$ such that
$\gs(0) = x(0)$ and $\gs (1) = y(0)$.
For each $n \ge 1$ let $l_n$ be the piecewise smooth loop 
$\gf^n (\gs) \cup -\gs$.

\begin{proposition}  \label{p:lin}
$\int_{l_n} \gl = nc$.
\end{proposition}

\proof
We start with
\begin{lemma}  \label{l:2}
$\int_{l_n} \gl = n \int_{l_1} \gl$.
\end{lemma}

\proof
Since $\gf$ is Hamiltonian, so is $\gf^k$ for any $k \ge 1$, and so
$\left\langle \Flux \gf^k, l_1 \right\rangle = 0$ for any $k \ge 1$.
Therefore,
\[
\int_{\gf^{k+1} \gs} \gl - \int_{\gf^k \gs} \gl 
\,=\, \int_{\gf^k l_1} \gl
\,=\, \int_{l_1} \left( \gf^k \right)^* \gl
\,=\, \int_{l_1} \gl 
\,=\, \int_{\gf \gs} \gl - \int_{\gs} \gl 
\]
and so
\begin{eqnarray*}
\int_{l_n} \gl \,=\, \int_{\gf^n \gs} \gl - \int_{\gs} \gl
&=& \sum_{k=0}^{n-1} \left( \int_{\gf^{k+1} \gs} \gl - \int_{\gf^k \gs} \gl
\right) \\ 
&=& \sum_{k=0}^{n-1} \left( \int_{\gf \gs} \gl - \int_{\gs} \gl \right)
\,=\, n \int_{l_1} \gl .
\end{eqnarray*}
\proofend

\begin{lemma}  \label{l:3}
$\int_{l_1} \gl = \ca_H(y) - \ca_H(x)$.
\end{lemma}

\proof
Define a $2$-chain $\Delta \colon [0,1] \times [0,1] \ra M$ by
$\Delta (t,s) = \gf_H^t \gs (s)$. 
Assuming the boundary of $[0,1] \times [0,1]$ to be oriented
counterclockwise we have
\[
\pp \Delta = - \gs - y + \gf \gs + x .
\]
A computation given in the proof of Proposition 2.4.A in \cite{P1} shows
that 
\[
\int_{\Delta} \omega 
\,=\, \int_0^1 H(t, x(t)) \, dt - \int_0^1 H(t, y(t)) \,dt .
\]
Putting everything together we conclude that
\begin{eqnarray*}
\int_{l_1} \gl
\,=\, \int_{\gf \gs} \gl - \int_{\gs} \gl
&=& \int_{\pp \Delta} \gl - \int_x \gl + \int_{y} \gl \\
&=& \int_{\Delta} \omega - \int_x \gl + \int_y \gl  
\,=\, \ca_H (y) - \ca_H (x) ,
\end{eqnarray*}
as claimed.
\proofend

We finally conclude from Lemmata \ref{l:2} and \ref{l:3} that 
\[
\int_{l_n} \gl \,=\, n \int_{l_1} \gl \,=\, nc ,
\]
and so the proof of Proposition \ref{p:lin} is complete.
\proofend

\m
\ni
We measure the lengths of curves in $M$ with respect to a 
Riemannian metric on $M$.
Let $x,y \in \cp_H$ and $l_n = \gf^n \gs \cup - \gs$ be as above. 
Choose $C < \infty$ so large that $\left| \gl (x) \right| \le C$ for all
$x \in \supp \gf \cup \gs$. Then $l_n \subset \supp \gf \cup
\gs$ for all $n \ge 1$ and so
\[
nc \,=\, \int_{l_n} \gl \,\le\, C \length l_n ,
\]
whence
$\frac{c}{C} n \le \length \gf^n \gs + \length \gs$ for all $n \ge 1$.
We conclude that
\begin{eqnarray*}
s_1(\gf) &\ge& \liminf_{n \ra \infty} \frac{1}{\log n} \log \length
               \gf^n \gs \\
         &\ge& \liminf_{n \ra \infty} \frac{1}{\log n} \log \left(
               \frac{c}{C} n \right) \\
         &=&   1.
\end{eqnarray*}

\b
\ni
{\bf Case 2: $\gf$ is not Hamiltonian.}
In this case, $\Flux \gf \in H_c^1(M;\RR)$ does not vanish.
Let $H_1^{cl}(M;\RR)$  be the first homology with closed support of $M$. Since
the pairing 
\[
H^1_c (M;\RR) \otimes H_1^{cl}(M;\RR) \ra \RR ,
\qquad 
\left( [\ga], [\gg] \right) \mapsto \int_{\gg} \ga ,
\]
is non-degenerate, we therefore find $[\gg] \in H_1^{cl} (M;\RR)$ 
such that 
\[
\left\langle \Flux \gf, [\gg] \right\rangle \,=\, \int_\gg \gf^* \gl - \gl
\,=:\, c \,>\, 0 .
\]
Since $\gf^* \gl - \gl$ has compact support, we can represent $[\gg]$ by
a smoothly embedded line $\gg \colon \RR \ra M$ such that $\gg \cap
\supp \gf \subset \gg \left( [0,1] \right)$. 
Define $\gs \colon [0,1] \ra M$ by $\gs (t) = \gg (t)$
and set again $l_n = \gf^n \gs \cup - \gs$.
Using that $\Flux$ is a homomorphism we then find 
\begin{eqnarray*}
\int_{l_n} \gl \,=\, \int_{\gs} \left(\gf^n\right)^* \gl - \gl 
               &=& \int_{\gg} \left(\gf^n\right)^* \gl - \gl \\
               &=& \left\langle \Flux \gf^n, [\gg] \right\rangle 
               \,=\, n \left\langle \Flux \gf, [\gg] \right\rangle 
               \,=\, n c .
\end{eqnarray*}
Proceeding as in Case 1 we conclude that $s_1 (\gf) \ge 1$.
The proof of Theorem~\ref{t:4exact} is complete.
\proofend

\subsection{Distortion in finitely generated subgroups of 
$\Symp_0^c \left(M,d\gl \right)$}
\label{ss:algebraic}
We consider an exact symplectic manifold $(M, d\gl)$ which is
convex at infinity. Theorem~\ref{t:4exact} yields at once
\begin{corollary}  \label{c:torsion}
The group $\Symp^c_0 (M, d\gl)$ has no torsion.
\end{corollary}

\proof
Assume that $\gf^m = \id$ for some $\gf \in \Symp^c_0 (M, d\gl)$
and $m \ge 1$. Using definition~\eqref{def:si} of $s_1(\gf)$ and
$\gf^{mn}(\gs) = 
\gs$ for all $n \ge 1$ and $\gs \in \Sigma_1$ we then find 
\begin{eqnarray*}
s_1 (\gf) &=&  
 \sup_{\gs \in \Sigma_1} \liminf_{n \ra \infty} 
    \frac{ \log \mu_g \left( \gf^n ( \gs ) \right)}{\log n} \\
          &\le&  \sup_{\gs \in \Sigma_1} \liminf_{n \ra \infty} 
    \frac{ \log \mu_g \left( \gf^{mn} ( \gs ) \right)}{\log mn} \\
          &=& 0,
\end{eqnarray*}
and so Theorem~\ref{t:4exact} implies $\gf = \id$.
\proofend

Proposition~\ref{p:4axa} can be used to obtain deeper insight into the
algebraic structure of the groups $\Ham^c \left(M,d\gl \right)$
and $\Symp_0^c \left(M,d\gl \right)$.
Following \cite{P1} we consider a finitely generated
subgroup $\cg$ of $\Symp^c_0 (M, d\gl)$.
Fix a set of generators of $\cg$ and denote by $\| \gf \|$
the word length of $\gf \in \cg$. The {\it distortion}\, $d(\gf)
\in [0,1]$ of $\gf \in \cg$ defined as
\[
d(\gf) \,=\, 1 - \liminf_{n \ra \infty} \frac{\log \| \gf^n
\|}{\log n }
\]
does not depend on the set of generators.

\begin{theorem}  \label{t:distortion}  
Consider a finitely generated subgroup $\cg$ 
of $\Symp^c_0 (M,d\gl)$.
\begin{itemize}
\item[(i)] 
If $\cg \subset \Ham^c (M, d\gl)$, then $d(\gf) =0$ for all
$\gf \in \cg \setminus \{ \id \}$. 
\item[(ii)] 
If $\cg \subset \Symp^c_0 (M, d\gl)$, then $d(\gf) \le \frac 12$
for all $\gf \in \cg \setminus \{ \id \}$, and $d(\gf)=0$ if
$\gf$ is not Hamiltonian.
\end{itemize}
\end{theorem}
\begin{question}  \label{q:12}
{\rm
Can the estimate $d(\gf) \le \frac 12$ in
Theorem~\ref{t:distortion}\:(ii) be replaced by $d(\gf) =0$?
}
\end{question}

\s \ni
{\it Sketch of the proof of Theorem~\ref{t:distortion}:}
The proof can be extracted from \cite{P1}, where an analogous
result was found for closed symplectic manifolds $(M, \go)$ with
$\pi_2(M)=0$. 
In our situation the arguments are considerably easier, however.
Choosing $k$ so large that $\cg \subset \Symp_0^c \left( M_k, d\gl \right)$
we can assume that $\left( M ,d\gl \right) = \left( M_k ,d\gl \right)$.
In addition to the geometric arguments in \cite[Sections
4.1--4.5]{P1} one uses that the action $\ca_H(x)$ of a
contractible $x \in \cp_H$ depends only on $x(0)$ and $\gf_H \in
\Ham^c \left( M_k, d\gl \right)$ 
(see \cite[Remark~3.1.1]{BPS} or \cite[Corollary~6.2]{FS}), 
that $\gl$ is bounded with respect to any Riemannian metric on
$M_k$, as well as Proposition~\ref{p:4axa},
Lemmata~\ref{l:2} and \ref{l:3}, and the flux homomorphism~\eqref{def:flux}.
\proofend

\ni
{\it Second proof of Corollary~\ref{c:torsion}:}
Let $\cg$ be the cyclic subgroup generated by 
$\gf \in \Symp^c_0 (M, d\gl)$.
If $\gf^m = \id$ for some $m \neq 0$, then $d(\gf)=1$, and so
Theorem~\ref{t:distortion} yields $\gf =\id$.
\proofend

Following again \cite{P1} we notice that
Theorem~\ref{t:distortion}\:(ii) can also be used to obtain
restrictions for representations of discrete groups on $\Symp^c_0
(M, d\gl)$.
An element $x$ of an abstract finitely generated group $\cg$ is
called a {\it
$\mathbf{U}$-element}\, if it is of infinite order and
\[
\liminf_{n \ra \infty} \frac{\log \| x^n \|}{\log n } \,=\,0 .
\]
Theorem~\ref{t:distortion}\:(ii) shows that $\phi (x) = \id$ for
every homomorphism $\phi \colon \cg \ra \Symp^c_0 (M, d\gl)$.
An example of a $\mathbf{U}$-element is the element $a$ of the
Baumslag--Solitar group
\[
BS (q,p) \,=\, \langle a,b \mid a^q = b a^p b^{-1} \rangle,
\quad\, q, p \in \ZZ \setminus \{ 0 \}, \; |p| < |q| ,
\]
see \cite[Example 1.6.E]{P1}.
Other examples of finitely generated groups containing a
$\mathbf{U}$-element are $SL (n;\ZZ)$ for $n \ge 3$ and, more generally,
irreducible non-uniform lattices in a semisimple real Lie group
whose real rank is at least two and which is connected, without
compact factors and with finite centre.
Let $\cg$ be such a lattice.
As in \cite[1.6.J]{P1} one obtains
\begin{corollary}  \label{c:finite}
Every homomorphism $\cg \ra \Symp_0^c (M, d\gl)$ has finite image.
\end{corollary}

\begin{remarks}\     
{\rm
{\bf 1.}
Combining the proof of Theorem~9.1.6 in \cite{MS3} with the
compactness theorems for $J$-holomorphic curves in geometrically
bounded symplectic manifolds \cite{G,Si} one sees that
Proposition~\ref{p:4axa} and hence Theorem~\ref{t:4exact} and the
results in Subsection~\ref{ss:algebraic} extend to geometrically bounded 
exact symplectic manifolds.

\s
\ni
{\bf 2.}
Consider the full group $\Symp^c \left( T^*S^d, d\gl \right)$ of
compactly supported symplectomorphisms of the standard cotangent bundle
$\left( T^*S^d, d \gl \right)$, where $d =1,2$. 
Recall that 
$\Symp^c \left( T^*S^d, d\gl \right) / \Symp_0^c \left( T^*S^d, d\gl \right)$
is infinite cyclic, \cite{S0}.
Corollary~\ref{c:torsion} thus implies that $\Symp^c \left( T^*S^d, d\gl
\right)$ has no torsion.
Moreover, given a lattice $\cg$ as in Corollary~\ref{c:finite}, 
every homomorphism $\cg \ra \ZZ$ is trivial because $\cg$ has
property ($T$), see \cite{DV}; together with Corollary~\ref{c:finite} we
thus see that every homomorphism $\cg \ra \Symp^c \left( T^*S^d, d \gl \right)$ has finite image.
}
\end{remarks}

\subsection{Proof of Proposition 1}  \label{prop1}

\ni
We first consider the open ball $B^{2n}(1)$ of radius $1$ in $\RR^{2n}$.
Choose a smooth function $f \colon \RR \ra [0,1]$ such that
\begin{equation*}  
 f(r) \,=\,
 \left\{
  \begin{array}{rll} 
     1   & \text{ iff} & r \le  \frac{1}{4n}, \\[0.2em]
     0   & \text{ iff} & r \ge 1 - \frac{1}{4n}, 
  \end{array}
 \right.
\end{equation*} 
and set $H(x) = f \left( \left| x \right|^2 \right)$.
Denoting by $J \in \cl \left( \RR^{2n} \right)$ the standard complex structure
on $\RR^{2n}$ we find that the time-$n$-map of the Hamiltonian flow of $H$
is given by
\[
\gf^n(x) \,=\, e^{2 f' \left( \left| x \right|^2 \right) n J} x, \quad \, x
\in B^{2n}(1) .
\]
Notice that $\gf^n$ preserves the Euclidean length of curves contained in a
round sphere. Moreover, there exists a constant $C < \infty$ such that
\[
\left\| d \gf^n (x) \right\| \le Cn  \quad \text{for all } x \in B^{2n}(1)
\text{ and } n \ge 1 ,
\] 
and so $s_i (\gf) \le 1$ for all $i$.

In order to show that $s_i(\gf) \ge 1$ for $i \in \left\{ 1, \dots, 2n-1
\right\}$ we shall investigate the $\gf$-orbits of particular cubes
$\gs_i$, which for convenience are defined on $Q_i' := \left[ 0,
\frac{1}{i+1} \right]^i$ instead on $[0,1]^i$. 
We measure the size of the cubes $\gf^n (\gs_i)$ with respect to the
Euclidean metric, and we use coordinates 
$\left( x_1, y_1, \dots, x_n,y_n \right)$ in $\RR^{2n}$.
We abbreviate 
\[
t = (t_1, \dots, t_i) \quad \text{ and } \quad
\bold{t} = \sqrt{t_1^2 + \dots +  t_i^2} .
\]

If $i \in \left\{ 1, \dots, 2n-1 \right\}$ is odd, we define $\gs_i
\colon Q_i' \ha B^{2n}(1)$ by
\begin{equation*}  
 \gs_i (t) \,=\,  
 \left\{
  \begin{array}{lll} 
     \left( t, 0, \dots, 0 \right)   
         & \text{ if} & i = 1, \\
     \left( t_1, 0, t_2, t_3, \dots, t_{i-1}, t_i, 0, \dots, 0 \right)
         & \text{ if} & i \ge 3 . 
  \end{array}
 \right.
\end{equation*} 
Then 
$\gf^n \left( \gs_i(t) \right) = e^{2f' \left(\bold{t}^2\right) n J} \gs_i(t)$.
A computation shows that
\[
\mu \left( \gf^n (\gs_i) \right) \,=\, 
\int_{Q_i'} \sqrt{ 1 + \big( 4 n f'' \left( \bold{t}^2 \right) t_1^2
\bold{t}^2 \big)^2 } \, dt .
\]
By our choice of $f$, this expression grows like $n$, and so $s_i (\gf) \ge 1$.

\s
If $i \in \left\{ 2, \dots, 2n-2 \right\}$ is even, we define $\gs_i
\colon Q_i' \ha B^{2n}(1)$ by
\begin{equation*}  
 \gs_i (t) \,=\,  
 \left\{
  \begin{array}{lll} 
     \left( t_1, 0, t_2, 0, \dots, 0 \right)   
         & \text{ if} & i = 2, \\
     \left( t_1, 0, t_2, 0, t_3, t_4, \dots, t_{i-1}, t_i, 0, \dots, 0 \right)
         & \text{ if} & i \ge 4 . 
  \end{array}
 \right.
\end{equation*} 
A computation shows that
\[
\mu \left( \gf^n (\gs_i) \right) \,=\, 
\int_{Q_i'} \sqrt{ 1 + \big( 4 n f'' \left( \bold{t}^2 \right) \left(
t_1^2 + t_2^2 \right) \bold{t}^2 \big)^2 } \, dt .
\]
By our choice of $f$, this expression grows like $n$, and so $s_i (\gf)
\ge 1$.

\s
Assume now that $(M, \go)$ is an arbitrary symplectic manifold.
By Darboux's theorem, there exists $\eps >0$ and a symplectic embedding
$\chi$ of the open ball $B^{2n}(\eps)$ of radius $\eps$ into $M$,
and proceeding as above we find a compactly supported Hamiltonian
diffeomorphism $\gf$ of $B^{2n}(\eps)$ such that $s_i(\gf) =1$ for all $i
\in \left\{ 1, \dots, 2n-1 \right\}$. 
Define the Hamiltonian diffeomorphism $\psi$ of $M$ by
\begin{equation*}  
 \psi (x) \,=\,
 \left\{
  \begin{array}{rll} 
     \chi \gf \chi^{-1}(x)    & \text{ if} & x \in \chi \left( B^{2n}
                                                     (\eps) \right), \\
     x   & \text{ if} & x \notin \chi \left( B^{2n} (\eps) \right) .
  \end{array}
 \right.
\end{equation*} 
Then $s_i(\psi) = s_i (\gf) =1$ for all $i \in \{ 1, \dots, 2n-1 \}$,
and so the proof of Proposition~1 is complete.
\proofend

\section{Proof of Theorem 2}  \label{t2}

\ni

\subsection{Unboundedness of the action functional}

Consider a closed manifold $B$ whose fundamental group
$\pi_1(B)$ is finite or has infinitely many conjugacy classes,
and let $H \colon S^1 \times T^*B \ra \RR$,
\[
H(t,q,p) \,=\, \tfrac12 \left| p-A(t,q) \right|^2 + V(t,q) ,
\]
be a classical Hamiltonian function.
Here, $\left| \;\! \cdot \;\! \right|$
refers again to the Riemannian metric on 
$T^*B$ induced by a fixed Riemannian metric on $B$.
The main ingredient of the proof of Theorem~2 is 

\begin{proposition}  \label{p:benci}
The action functional $\ca_H$ is not bounded from above on $\cp_H$.
\end{proposition} 

\proof
Since $H$ is fibrewise convex, its Legendre transform $L \colon
S^1 \times TB \ra \RR$ is defined and equals
\begin{equation}  \label{def:Lag}
L(t,q, \dot{q}) \,=\, \tfrac 12 \left| \dot{q} \right|^2 +
\langle B(t,q), \dot{q} \rangle + W(t,q)
\end{equation}
where $B(t,q) = -A(t,q)$ and $W(t,q) = \frac 12 \left| A(t,q)
\right|^2 - V(t,q)$.
The $1$-periodic orbits $x(t)$ of the flow of $H$ correspond to
$1$-periodic orbits $q(t)$ of the Lagrangian flow generated by
$L$,
and $\ca_H(x)$ equals the Lagrangian action
\[
\cl_H (q) \,=\, \int_0^1 L \left( t,q,\dot{q} \right) dt.
\]
It is proved in \cite{Be}\footnote{Lemma~2.6\,(c) in \cite{Be}
should, however, read
\[
\int_0^1 \ce \left( \gb (\gl) \right)^{1/2} d\gl \,\le\, d_\gb +
\ce \left( \gb (0) \right)^{1/2} ,
\]
and then the proofs of Lemmata~2.7 and 4.3 should be corrected accordingly.} that if
$\pi_1(B)$ is finite, then $\cl_H$ is not bounded from above on
the set of $1$-periodic orbits, and so the same holds for
$\ca_H$.
Assume now that $\pi_1(B)$ has infinitely many conjugacy
classes.
These conjugacy classes correspond to the connected components
of the space $\Omega B$ of continuous maps $S^1 \ra B$, and
$\Omega B$ is homotopy equivalent to the Sobolev space
$\Omega^1B = W^1 \left( S^1, B \right)$ of maps $S^1 \ra B$ with
``square integrable derivative'', see \cite{Be}.
It is proved in \cite{Be} that $\cl_H$ is a
$C^1$-functional on $\Omega^1B$ and that each of the infinitely
many components of
$\Omega^1B$ contains a critical point $q_n$ of $\cl_H$, $n=1,2,\dots$.
Since these loops are continuous, they are smooth $1$-periodic
orbits of the Lagrangian flow.
Arguing by contradiction we assume that $\left\{ \cl_H (q_n)
\right\} \subset \RR$ is bounded from above.
According to Lemma~2.4 in \cite{Be}, $\cl_H$ is bounded from
below on $\Omega^1B$, and so $\left\{ \cl_H (q_n) \right\}$ is
bounded;
moreover, $\nabla \cl_H (q_n) =0$ for all $n$.
Since $\cl_H$ satisfies the Palais--Smale condition on $\Omega^1
B$, \cite{Be}, we
conclude that the sequence $(q_n)$ has a convergent subsequence
in $\Omega^1B$, which is impossible.
This contradiction completes the proof of Proposition~\ref{p:benci}.
\proofend

\subsection{End of the proof}
Let $B$ and $H$ be as in Theorem~2.
In view of Proposition~\ref{p:benci} we find $x,y \in \cp_H$ such that $c:=
\ca_H(y)-\ca_H(x) >0$.
Choose a smoothly embedded curve $\gs \colon [0,1] \ra T^*B$ such that
$\gs(0) = x(0)$ and $\gs (1) = y(0)$.
For each $n \ge 1$ let $l_n$ be the piecewise smooth loop 
$\gf_H^n (\gs) \cup -\gs$. According to Proposition~\ref{p:lin} we then have
\begin{equation}  \label{e:lnnc}
\int_{l_n} \gl \,=\, nc .
\end{equation}
We measure the lengths of curves in $T^*B$ with respect to $g^*$. 

\m
\ni
{\it Proof of Theorem~2 (i).} 
We choose $R>0$ so large that
\begin{equation}  \label{inc:xy}
x(0), y(0) \in T_R^* B.
\end{equation}
\begin{claim}  \label{claim}
$\length l_n \ge \min (c,1) \sqrt{n}$ for any $n \ge 4 R^2$.
\end{claim}
\proof
Fix $n \ge 4R^2$.
We distinguish two cases.

\m
\ni
Case 1: $l_n \subset T^*_{\sqrt{n}} B$.

\s
\ni
We have $\left| \gl \left( l_n(t) \right) \right| \le \sqrt{n}$ for all $t
\in [0,1]$. Together with the identity~\eqref{e:lnnc} we infer that
\[
nc \,=\, \int_{l_n} \gl \,\le\, \sqrt{n} \length l_n ,
\]
and so $\length l_n \ge c \sqrt{n}$.

\m
\ni
Case 2: $l_n(t) \notin T^*_{\sqrt{n}} B$ for some $t \in [0,1]$.

\s
\ni
In view of \eqref{inc:xy} and the assumption $n \ge 4R^2$ we can
estimate 
\begin{eqnarray*}
\length l_n  &\ge& \left( \left| p \left( l_n(t) \right) \right| - 
                            \left| p (x(0)) \right| \right) +
                     \left( \left| p \left( l_n(t) \right) \right| - 
                            \left| p (y(0)) \right| \right)  \\
             &\ge& \left( \sqrt{n} -R \right) +  \left( \sqrt{n} -R
                                                     \right) \\
             &\ge& \frac{\sqrt{n}}{2} + \frac{\sqrt{n}}{2} ,
\end{eqnarray*}
and so $\length l_n \ge \sqrt{n}$.
\proofend

According to Claim \ref{claim}, 
\[
\min (c,1) \sqrt{n} \,\le\, \length \gs + \length \gf_H^n \gs
\]
for all $n \ge 4R^2$, and so
\begin{eqnarray*}
s_1(\gf) &\ge& \liminf_{n \ra \infty} \frac{1}{\log n} \log \length
               \gf_H^n \gs \\
         &\ge& \lim_{n \ra \infty} \frac{1}{\log n} \log \left(
               \min (c,1) \sqrt{n} \right) \\
         &=&   \frac{1}{2}.
\end{eqnarray*}

\m
\ni
{\it Proof of Theorem~2 (ii).} 
Let $r < \infty$ be such that $H$ is
independent of $t$ on $T^*B \setminus T^*_r B$.
Choosing $r$ larger if necessary we can assume that $\gs \subset T^*_r B$.
Choose $a$ so large that $T_r^*B$ is contained in the sublevel set 
\[
H^a \,=\, \left\{ (t, q,p) \in S^1 \times T^*B \mid H(t,q,p) \le a
\right\} .
\]
Since $H$ is time-independent on the boundary of $H^a$, the set $H^a$ is
invariant under $\gf_H$, and so $\gf_H^n \gs \subset H^a$ for all $n$.
Choosing $R < \infty$ so large that $H^a \subset T_R^*B$, we then have $l_n
\subset T_R^*B$ for all $n \ge 1$.
In particular, $\left| \gl \left( l_n(t) \right) \right| \le R$ 
for all $n \ge 1$ and all $t \in [0,1]$. 
Together with \eqref{e:lnnc} we infer that
$\length l_n \ge \frac{c}{R} n$ for all $n \ge 1$.
Proceeding as in the proof of (i) we conclude that $s_1(\gf_H) \ge 1$.
The proof of Theorem~2 is complete. 
\proofend

\begin{remarks}\  \label{r:general}
{\rm
(i)
The results from \cite{Be} used in the proof of
Proposition~\ref{p:benci}
hold for more general convex Lagrangians than classical ones of the form
\eqref{def:Lag}, and so Proposition~\ref{p:benci} and hence
Theorem~2 hold for more general convex Hamiltonians.

\s
(ii)
Using \cite{C} and \cite[Section~6]{CS1} supplemented by a recent result in
\cite{AS}, one finds that Proposition~\ref{p:benci} and hence
Theorem~2 hold for a yet larger class of Hamiltonians $H$, which do
not need to be convex.
They are only assumed to satisfy the two asymptotic conditions
\begin{itemize}
\item[(H1)]
$dH (X) (t,q,p) - H(t,q,p) \,\ge\, c |p|^2 - C$,
\vspace{0.3em}
\item[(H2)]
$\left| \nabla_q H(t,q,p) \right| \le 
C \left( 1+ |p|^2 \right)$
\; and \;
$\left| \nabla_p H(t,q,p) \right| \le 
C \left( 1+ |p| \right)$,
\end{itemize}
for all $(t,q,p) \in S^1 \times T^*B$. 
Here, $X= \sum_i p_i \frac{\pp}{\pp p_i}$ is the Liouville vector field,
$\nabla$ denotes the Levi--Civita connection with respect
to the Riemannian metric $g$ on $B$, and $c$ and $C$ are positive
constants depending only on $g$.
}
\end{remarks}

\section{Proof of Theorem 3}  \label{t3}

\ni
In this section we shall first describe the known Riemannian manifolds
with periodic geodesic flow
and then define twists on such manifolds.
We then prove Theorem~3 and Corollary~3, and finally study twists from a
topological point of view.

\subsection{$P$-manifolds}  \label{Pmanifolds}

\ni
Geodesics of a Riemannian manifold will always be parametrized by arc-length.
A {\it $P$-manifold}\, is by definition a connected Riemannian
manifold all whose geodesics are periodic. It follows from Wadsley's
Theorem that the 
geodesics of a $P$-manifold admit a common period, 
see \cite{Wy} and \cite[Lemma 7.11]{B}.
We normalize the Riemannian metric such that the minimal common period is $1$. 
Every $P$-manifold is closed, and besides $S^1$ every $P$-manifold 
has finite fundamental group, see \cite[7.37]{B}.
The main examples of $P$-manifolds are the CROSSes
\[
S^d, \quad \RP^d, \quad \CP^n, \quad \HP^n, \quad \Ca 
\]
with their canonical Riemannian metrics suitably normalized. 
The simplest way of obtaining other $P$-manifolds is to look at
Riemannian quotients of CROSSes. The main examples thus obtained are
the spherical space forms $S^{2n+1} / G$ where $G$ is a finite
subgroup of $O (2n+2)$ acting freely on $S^{2n+1}$.
These spaces are classified in \cite{W}, and examples are lens spaces,
which correspond to cyclic $G$.
According to \cite[pp.\ 11--12]{AM} and \cite[7.17\, (c)]{B}, 
the only other Riemannian quotients of CROSSes are the spaces
$\CP^{2n-1} / \ZZ_2$; here, the fixed point free involution on
$\CP^{2n-1}$ is induced by the involution
\begin{equation}  \label{e:inv2}
(z_1, z_1', \dots, z_n, z_n') \,\mapsto\,
\left( \bar{z}'_1, -\bar{z}_1, \dots, \bar{z}'_n, -\bar{z}_n \right)
\end{equation}
of $\CC^{2n}$. Notice that $\CP^1 / \ZZ_2 = \RP^2$. We shall thus assume
$n \ge 2$.
On spheres, there exist $P$-metrics which are not isometric to the round
metric $g_{\can}$.
We say that a $P$-metric on $S^d$ is a {\it Zoll metric}\, if it can be
joined with $g_{\can}$ by a smooth path of $P$-metrics.
All known $P$-metrics on $S^d$ are Zoll metrics.
For each $d \ge 2$, the Zoll metrics on $S^d$ form an infinite
dimensional space. 
For $d \ge 3$, the known Zoll metrics admit $SO(d)$ as isometry group,
but for $d=2$, the set of Zoll metrics contains an open set all of whose
elements have trivial isometry group.
We refer to \cite[Chapter 4]{B} for more information about Zoll metrics.

CROSSes, their quotients and Zoll manifolds are the only known
$P$-manifolds.
It would be interesting to know whether this list is complete.
As an aside, we mention that for the known $P$-manifolds all
geodesics are simply closed. Whether this is so for all
$P$-manifolds is unknown, \cite[7.73 (f''\!')]{B}, except for
$P$-metrics on the $2$-sphere, \cite{GG}.

\m
For a geodesic $\gg \colon \RR \ra B$ of a $P$-manifold $(B,g)$ and $t>0$
we let $\ind \gg (t)$ be the number of linearly independent Jacobi
fields along $\gg (s)$, $s \in [0,t]$, which vanish at $\gg (0)$ and
$\gg (t)$.
If $\ind \gg (t) >0$, then $\gg (t)$ is said to be conjugate to $\gg
(0)$ along $\gg$.
The index of $\gg$ defined as
\[
\ind \gg \,=\, \sum_{t \in ]0,1[} \ind \gg (t) 
\]
is a finite number.
According to \cite[1.98 and 7.25]{B}, every geodesic on $(B,g)$ has the
same index, say $k$.
We then call $(B,g)$ a {\it $P_k$-manifold}.
\begin{proposition}  \label{p:index}
For the known $P$-manifolds, the indices of geodesics are as follows.
\begin{table}[h]  \label{ta:ind}
 \begin{center}
 \renewcommand{\arraystretch}{1.5}
  \begin{tabular}{|c||c|c|c|c|c|} \hline
  $(B,g)$ & $S^d$ & $\RP^d$ & $\CP^n$ & $\HP^n$ & $\Ca$   \\ \hline
  $k$     & $d-1$ & $0$     & $1$     & $3$     & $7$     \\ \hline
  \end{tabular}
 \end{center}
\end{table} 

\ni
For a quotient $S^{2n+1} / G$ we have $k=0$ if $- \idd \in G$ and $k =
2n$ otherwise,
and for $\CP^{2n-1} / \ZZ_2$ we have $k = 1$.
Finally, for a Zoll manifold modelled on $S^d$ we have $k=d-1$.
\end{proposition} 

\proof
For CROSSes, the result is well known, see \cite[3.35 and 3.70]{B}.
For this proof, we assume the Riemannian metrics on $S^{2n+1}/ G$ and
$\CP^{2n-1}/\ZZ_2$ to be locally isometric to the normalized Riemannian
metrics on $S^{2n+1}$ and $\CP^{2n-1}$, respectively.
The Jacobi fields on $S^{2n+1}/ G$ and $\CP^{2n-1}/\ZZ_2$ then
correspond to Jacobi fields on $S^{2n+1}$ and $\CP^{2n-1}$. It follows
that if all geodesics on $S^{2n+1}/ G$ of length $1/2$ are closed, i.e.,
$- \id \in G$, then $k=0$, and $k=2n$ otherwise.
The isometry $\gs$ of $\CP^{2n-1}$ induced by \eqref{e:inv2} maps a
point $x$ to its conjugate locus diffeomorphic to $\CP^{2n-2}$, and
the family of geodesics emanating from $x$ which pass through $\gs(x)$ 
is $2$-dimensional.
Since $n \ge 2$, there are geodesics not passing through $\gs(x)$, and
so $k=1$. 
Finally, let $g$ be a Zoll metric on $S^d$. In order to show $k=d-1$, we
choose a smooth family $g_t$ of $P$-metrics on $S^d$ such that 
$g_0 = g_{\can}$ and $g_1 = g$.
Fix $x \in S^d$ and $v \in T_xS^d$ with $|v|_0=1$.
For $t \in [0,1]$ let $\gg_t$ be the geodesic of 
$( S^d, g_t )$ with $\gg_t (0)=x$ and $\dot{\gg}_t (0) =
\frac{|v|_0}{|v|_t}v$.
The set
$T = \left\{ t \in [0,1] \mid \ind \gg_t = d-1 \right\}$
contains $0$, and comparing Jacobi fields on $( S^d, g_t )$
along $\gg_t$ one easily shows that $T$ is both open and closed.
In particular, $\ind \gg_1 = d-1$, and so $( S^d, g_1 )$
is a $P_{d-1}$-metric.
\proofend

We conclude this subsection by proving a property shared by all
$P$-manifolds $(B,g)$.
For each point $v$ in the unit tangent bundle $\pp T_1B$ over $B$ 
we denote by $l(v)$ the length of the simply closed orbit of the 
geodesic flow through $v$.
By our normalization of $g$, $l(v) = 1/j$ for some integer $j$ depending on $v$.

\begin{lemma}  \label{l:dicht}
The set of point $v \in \pp T_1B$ with $l(v) =1$ is open and dense in 
$\pp T_1B$.
\end{lemma}

\proof
According to \cite[Corollary A.18]{B}, the length function $l$ is
continuous on an open and dense subset $V$ of $\pp T_1B$.
It follows from Proposition~4.5 of \cite{W} that $l |_V \equiv 1/j$ 
for some integer $j$
and that for each $v \in \pp T_1B$ there is an integer $h$ such that $l(v) =
1/ (hj)$.
Our choice of $g$ implies that $j=1$, and so $l(v)=1$ for all $v \in V$.
\proofend

\subsection{Twists}  \label{twists}

\ni
Consider a $P$-manifold $(B,g)$.
As before, we choose coordinates $(q,p)$ on $T^*B$, and using $g$ we
identify the cotangent bundle $T^*B$ with the tangent bundle $TB$. The
Hamiltonian flow of the function $\frac{1}{2} |p|^2$ corresponds to the
geodesic flow on $TB$.
For any smooth function $f \colon [0,\infty[\; \ra [0,\infty[$ such that 
\begin{equation}  \label{e:f}
\text{$f(r)=0$ for $r$ near $0$ \; and \; $f'(r)=1$ for $r \ge 1$}
\end{equation}
we define the {\it twist}\, $\gt_f$ as the time-$1$-map of the Hamiltonian
flow generated by $f \left( |p| \right)$.
Since $(B,g)$ is a $P$-manifold, $\gt_f$ is the identity on $T^*B
\setminus T_1^*B$, and so $\gt_f \in \Symp^c \left( T^*B \right)$.

\begin{proposition}  \label{prop:n}
(i) The class $[\gt_f] \in \pi_0 \left( \Sympcc \left(T^*B \right) \right)$ 
does not depend on the choice of $f$.

\s
(ii)
$s_i \left( \gt_f^m \right) = l \left( \gt_f^m \right) = 1$ 
for every $i \in \left\{ 1, \dots, 2d-1 \right\}$, 
every $m \in \ZZ \setminus \{0\}$ and every $f$.
\end{proposition}

\proof
(i)\,
Let $f_i \colon [0,1] \ra [0,1]$, $i=1,2$, be two functions as in \eqref{e:f}. 
Then the functions $f_s = (1-s) f_1 + s f_2$, $s \in [0,1]$, 
are also of this form, and $s \mapsto \gt_{f_s}$ is an isotopy in
$\Sympcc \left( T^*B \right)$ joining $\gt_{f_1}$ with $\gt_{f_2}$.

\s
(ii)\,
The proof is similar to the proof of Proposition~1.
Without loss of generality we assume $m=1$.
Let $\gt^t$ be the Hamiltonian flow of $f \left( |p| \right)$.
Then $\gt_f = \gt^1$.
For each $r>0$ the hypersurface $S_r = \pp T_r^*B$ is invariant
under $\gt^t$.
We denote by $\gt_r^t$ the restriction of $\gt^t$ to $S_r$.
As before, we endow $T^*B$ with the Riemannian metric $g^*$.
For $x \in S_r$ let $\left\| d \gt^t_r (x) \right\|$ be the operator
norm of the differential of $\gt_r^t$ at $x$ induced by $g^*$.
Since $\gt_1^t$ is $1$-periodic, we find $C<\infty$ such that 
$\left\| d \gt^t_1(x) \right\| \le C$ for all $t$ and all $x \in S_1$.
Since
\[
\gt_r^t (x) \,=\, r \gt_1^{f'(r)t} \left( \tfrac{x}{r} \right)
\]
for all $t \in \RR$, $r>0$ and $x \in S_r$, we conclude that
\begin{equation}  \label{e:t3:C}
\left\| d \gt_r^t (x) \right\| \,=\, 
\left\| d \gt_1^{f'(r)t} \left( \tfrac{x}{r} \right) \right\| \,\le\, C
\end{equation}
for all $t \in \RR$, $r>0$ and $x \in S_r$.
We next fix $(q,p) \in T^*B \setminus B$
and consider the line $F_p = \RR p \subset T_q^*B$ orthogonal to
$S_{|p|}$ through $(q,p)$.
We denote by $\gt_p^t$ the restriction of $\gt^t$ to $F_p$.
Let $\gg$ be the geodesic on $B$ with $\gg (0) = q$ and $\dot{\gg} (0) =
  \frac{p}{|p|}$.
Then $\gt_p^t(F_p) \subset T^*\gg$ for all $t$. 
As a parametrized curve, $\gg$ is isometric to the circle $S^1$
of length $1$, and so $T^*\gg \setminus \gg$ is isometric to
$T^*S^1 \setminus S^1$. We thus find 
\begin{equation}  \label{e:t3:C'}
\left\| d \gt_p^t(q,p) \right\| \,=\,
\sqrt{ \left( f''(|p|) t \right)^2 +1 } \,\le\, C' (t+1) 
\end{equation}
for all $t \ge 0$ and $(q,p) \in T^*B$ and some constant $C' < \infty$.
The estimates \eqref{e:t3:C} and \eqref{e:t3:C'} show that for any
$i$-cube $\gs \colon Q^i \ha T^*B$, 
\[
\mu_{g^*} \left( \gt_f^n (\gs) \right) \,\le\, C^{i-1} C' (n+1)
\]
for all $n \ge 1$, and so $s_i (\gt_f) \le 1$ for all $i$.

\s
We are left with showing $s_i(\gt_f) \ge 1$ for 
$i \in \left\{ 1, \dots, 2d-1 \right\}$.
Fix $q \in B$, choose an orthonormal basis $\left\{ e_1, \dots, e_d \right\}$ of
$T_q^*B$, and let $R$ be such that $f \left( [0,1] \right) \subset
[0,R]$.
For $j \in \left\{1, \dots, d \right\}$ we let $E^j$ be the subspace
of $T^*_qB$ generated by $\left\{ e_1, \dots, e_j \right\}$, 
and we set $E_R^j = E^j \cap T^*_RB$.
We first assume $i \in \left\{ 1, \dots, d \right\}$. 
We then choose $\gs \colon Q^i \ha E^i$ such that $E_R^i \subset
\gs (Q^i)$.
The set $\pi \left( \gt_f (\gs) \right)$ consists of a smooth
$(i-1)$-dimensional family of geodesics in $B$. Its $i$-dimensional
measure $\mu_g \left( \pi \left( \gt_f (\gs) \right) \right)$ with
respect to $g$ thus exists and is positive.
Moreover, 
$\pi \left( \gt_f^n (\gs) \right) = \pi \left( \gt_f (\gs) \right)$
for all $n \ge 1$, and every point in $\pi \left( \gt_f^n (\gs)
\right)$ has at least $2n$ preimages in $\gt_f^n (\gs)$.
Since $\pi \colon \left( T^*B, g^* \right) \ra (B,g)$ is a Riemannian
submersion, we conclude that
\[
\mu_{g^*} \left( \gt_f^n (\gs) \right) \,\ge\, 
2n \mu_g \left( \pi \left( \gt_f^n (\gs) \right) \right) ,
\]
and so $s_i (\gt_f) \ge 1$.

Since the fibre $T_q^*B$ is a Lagrangian submanifold of $T^*B$, we have
shown that $s(\gt_f) = l(\gt_f)=1$.
We shall therefore only sketch the proof of the remaining inequalities
$s_i (\gt_f) \ge 1$, $i \in \left\{ d+1, \dots, 2d-1 \right\}$.
For such an $i$ we choose a small $\eps >0$, set $B_i = \exp_q
E_{\eps}^{i-d}$, and choose $\gs \colon Q^i \ha T^*B_i$ such that
$T_R^*B_i \subset \gs \left(Q^i \right)$.
Then there is a constant $c>0$ such that 
\begin{equation}  \label{est:mgt}
\mu_{g^*} \left( \gt^n_f (\gs) \right) \,\ge\, cn 
            \quad \, \text{for all } n \ge 1 .
\end{equation}
This is so because $\gt_f$ restricts to a symplectomorphism on the $i-d$
cylinders $T^*\gg_j$ over the geodesics $\gg_j$ with $\gg_j(0)=q$ and
$\dot{\gg} (0) = e_j$, and -- as we have seen above -- grows linearly on
the $(2d-i)$-dimensional remaining factor in the fibre.
An explicit proof of \eqref{est:mgt} can be given by computing the
differential $d \gt_f(q,p)$ with respect to suitable orthonormal bases
of $T_{(q,p)}T^*B$ and $T_{\gt_f(q,p)} T^*B$.
\proofend

Let $\tau \in \Symp^c \left( T^* S^d \right)$ be a generalized Dehn
twist as defined in Figure~\ref{figure2}; for an analytic definition we
refer to \cite[5a]{S2}.
Then $\tau^2$ is a twist $\gt_f$.
Proposition~\ref{prop:n} (ii) and the argument given in \ref{c:3} below thus
show that
$s_i \left( \tau^m \right) = l \left( \tau^m \right) = 1$ 
for every $i \in \left\{ 1, \dots, 2d-1 \right\}$ and 
every $m \in \ZZ \setminus \{0\}$.

\b
Theorem~3 is a special case of the following theorem, which is the main
result of this section.
\begin{theorem}  \label{t:twist}
Let $(B,g)$ be a $d$-dimensional $P_k$-manifold, and let $\gt$ be a
twist on $T^*B$.
If $d=2$ and $B$ is diffeomorphic to $\RP^2$, assume that $g =
g_{\can}$, and if $d \ge 3$ and $k=1$, assume that $(B,g)$ is 
$\CP^n$ or $\CP^{2n-1} / \ZZ_2$.
If $\gf \in \Symp^c \left(T^*B \right)$ is such that 
$[ \gf ] = \left[ \gt^m \right] \in \pi_0 \left( \Symp^c \left( T^*B
\right) \right)$ for some $m \in \ZZ \setminus \{0\}$, 
then $s_d (\gf ) \ge l(\gf) \ge 1$.
\end{theorem} 

\begin{remarks}\
{\rm 
{\bf 1.}
Theorem~\ref{t:twist} covers all known $P$-manifolds.
Indeed, the only known $P$-metric on $\RP^2$ is $g_{\can}$, 
and the only known $P_1$-manifolds of dimension at least $3$ are 
$\CP^n$ and $\CP^{2n-1} / \ZZ_2$.
The following two results suggest that Theorem~\ref{t:twist}
in fact covers all $P$-manifolds.

\s
(i)
If $g$ is a $P$-metric on $\RP^2$ such that for some point $x$ there
exists $l>0$ such that all geodesics of length $l$ emanating from $x$
are embedded circles, then $g = g_{\can}$ by Green's theorem \cite{Gr}.

\s
(ii)
If $(B,g)$ is a $P_1$-manifold containing a point $x$ for which
there exists $l>0$ such that for each geodesic $\gg$ emanating from
$x$,
$\gg (l) = x$ and $\gg (t) \neq x$ for all $t \in (0,l)$,
then $B$ has the homotopy type of $\CP^n$, see \cite[7.23]{B}.

\s
{\bf 2.}
(i)
We recall from \cite{S0,S8} that 
$\pi_0 \left( \Symp^c \left( T^*S^2 \right) \right) = \ZZ$ 
is generated by the class $[\tau]$ of a generalized Dehn twist.

\s
(ii)
For $S^d$, $d \ge 3$, $[\tau]^2 = [\gt]$ and Theorem~\ref{t:twist} imply
that $[\tau]$ generates an infinite cyclic subgroup of 
$\pi_0 \left( \Symp^c \left( T^*S^d \right) \right)$,
and for those $P_k$-manifold $(B,g)$ covered by
Theorem~\ref{t:twist} which are not diffeomorphic to a sphere, 
Theorem~\ref{t:twist} implies that
$[\gt]$ generates an infinite cyclic subgroup of
$\pi_0 \left( \Symp^c \left( T^*B \right) \right)$.
This was proved in \cite[Corollary~4.5]{S2} {\it for all}\,
$P$-manifolds.\footnote{It is assumed in \cite{S2} that
$H^1(B;\RR)=0$. Besides for $B=S^1$ this is, however, guaranteed
by the Bott--Samelson Theorem \cite[Theorem~7.37]{B}.}
It would be interesting to know whether there are other elements in
these symplectic mapping class groups.
\diam
}
\end{remarks}

Theorem~\ref{t:twist} is proved in the next two subsections.

\subsection{Lagrangian Floer homology}

\ni
Floer homology for Lagrangian submanifolds was invented by Floer
in a series of seminal papers, \cite{F2,F3,F1,F4},
and more general versions have been developed meanwhile,
\cite{Oh1, FOOO}.
In this subsection we first follow \cite{KS} and define
Lagrangian Floer homology for certain pairs of 
Lagrangian submanifolds with
boundary in an exact compact convex symplectic manifold.
We then compute this Floer homology in the special case that the
pair consists of
a fibre and the image of another fibre
under an iterated twist on the unit coball bundle over a $P$-manifold.

\subsubsection{Lagrangian Floer homology on convex symplectic manifolds}
\label{ss:LagFloer}
We consider an exact compact connected symplectic manifold 
$\left( M , \go \right)$ with boundary $\pp M$ and two compact Lagrangian
submanifolds $L_0$ and $L_1$ of $M$ meeting the following hypotheses.
\begin{itemize}
\item[\bf{(H1)}\,] $L_0$ and $L_1$ intersect transversally; 
\item[\bf{(H2)}\,] $L_0 \cap L_1 \cap \partial M = \emptyset$;
\item[\bf{(H3)}\,] $H^1 (L_j; \RR) = 0$\, for $j=0,1$.
\end{itemize}
We also assume that there exists a Liouville vector field $X$
(i.e., $\cl_X \go = d \iota_X \go = \go$) which is defined on a
neighbourhood $U$ of $\pp M$ and is everywhere transverse to
$\pp M$, pointing outwards, such that
\begin{itemize}
\item[\bf{(H4)}\,]
 $X(x) \in T_x L_j \quad \text{ for all }\, x \in L_j \cap U,\; j = 0,1.$
\end{itemize}
Let $\gf_r$ be the local semiflow of $X$ defined near $\pp M$.
Since $\pp M$ is compact, we find $\eps >0$ such that $\gf_r(x)$
is defined for $x \in \pp M$ and $r \in [-\eps,0]$. For these
$r$ we set
\[
U_r \,=\, \bigcup_{r \le r' \le 0} \gf_{r'} \left( \pp M \right)
.
\]
In view of (H2) there exists 
$\epsilon' \in \left] 0,\epsilon \right[$ such that for 
$V = U_{\eps'}$ we have
\begin{equation}  \label{e:empty}
V \cap L_0 \cap L_1 = \emptyset .
\end{equation}
An almost complex structure $J$ on $(M,\go)$ is called {\it
$\go$-compatible}\, if 
$\go \circ \left( \id \times J \right)$ is a Riemannian metric
on $M$.
Following \cite{CFH,V2,BPS} we consider the space $\cj$ of
smooth families $\bJ = \{ J_t\}$, $t \in [0,1]$, 
of smooth $\go$-compatible almost complex structures on $M$
such that $J_t(x) = J(x)$ does not depend on $t$ for $x \in V$
and such that
\begin{itemize}
 \item[(J1)] 
 $\go \left( X(x), J(x) v \right) = 0,
 \qquad 
 x \in \pp M, \,\, v \in T_x \pp M$, \vspace{0.4em}
 \item[(J2)] 
 $\go \left( X(x), J(x) X(x) \right) = 1, 
 \qquad 
 x \in \pp M$, \vspace{0.4em}
 \item[(J3)] 
 $d\gf_r(x) J(x) = J (\gf_r(x)) d \gf_r(x), 
 \qquad 
 x \in \pp M, \,\, r \in [-\eps',0]$.
\end{itemize}
For later use we examine conditions (J1) and (J2) closer.
The contact structure $\xi$ on $\pp M$ is defined as
\begin{equation}  \label{def:xi}
\xi \,=\, \left\{ v \in T\pp M \mid \go (X,v) =0 \right\} ,
\end{equation}
and the Reeb vector field $R$ on $\pp M$ is defined by
\begin{equation}  \label{def:Reeb}
\go (X,R) =1
\qquad \text{ and } \qquad
\go (R,v) =0 \;\text{ for all }\, v \in T\pp M .
\end{equation}
\begin{lemma}  \label{l:equiv}
Conditions (J1) and (J2) are equivalent to 
\[
J \xi = \xi 
\qquad \text{ and } \qquad
JX=R .
\]
\end{lemma}
\ni
The proof follows from definitions and the $J$-invariance of $\go$.
It follows from Lemma~\ref{l:equiv} that the set $\cj$ is
nonempty and connected, see \cite{CFH}.
Let
\[
S \,=\, \left\{ z=s+it \in \CC \mid s \in \RR, t \in [0,1] \right\} 
\]
be the strip. 
The energy of $u \in C^\infty \left( S,M \right)$ is defined as
\[
E(u) \,=\, \int_{S} u^* \go .
\]
For $u \in C^\infty \left( S,M \right)$ consider Floer's 
equation
\begin{eqnarray}  \label{floer}
 \left\{       
  \begin{array}{lcr}
    \pp_s u+J_t(u) \pp_t u = 0,  \\ [0.2em]
    u(s,j) \in L_j,\,\, j \in \{0,1\}, \\ [0.2em]
    E(u) < \infty .
  \end{array}
 \right.
\end{eqnarray}
Notice that for a solution $u$ of \eqref{floer},
\[
E(u) \,=\, \int_{S} \left\| \pp_s u \right\|^2 \,=\, \frac 12 \int_S
\left\| \pp_s u \right\|^2 + \left\| \pp_t u \right\|^2 
\]
is the energy of $u$ associated with respect to any Riemannian metric
defined via an $\go$-compatible $J$.
It follows from (H1) that for every solution $u$ of 
\eqref{floer} there exist
points $c_-, c_+ \in L_0 \cap L_1$ such that 
$\lim_{s \ra \pm \infty} u(s,t) = c_{\pm}$ uniformly in $t$, 
cf.\ \cite[Proposition~1.21]{Sa}.
The following lemma taken from \cite{EHS,KS} shows that 
the images of solutions of \eqref{floer} uniformly stay away from 
$\pp M$.
\begin{lemma}  \label{l:stayaway}
Let $u$ be a finite energy solution of \eqref{floer}. Then 
\[
u(S) \cap V = \emptyset .
\]
\end{lemma}

\proof
Define $f \colon V \ra \RR$ by $f \left( \gf_r(x) \right) =
e^r$,
where $x \in \pp M$ and $r \in [-\eps',0]$.
Using (J1), (J2), (J3) we find that
the gradient $\nabla f$ with respect to each metric
$\go \circ \left( \id \times J_t \right)$ is $X$;
for the function
\[
F \colon \Omega = u^{-1}(V) \to \RR, \quad\, (s,t) = z \mapsto F(z) = f
\circ u(z) ,
\]
one therefore computes $\Delta F = \langle \pp_s u, \pp_s u
\rangle$, see e.g.\ \cite{FS}, so that $F$ is subharmonic.
It follows that $F$ does not attain a strict maximum on
the interior of $\Omega$.
In order to see that this holds on $\Omega$, fix a point $z
\in \pp S$. We first assume $z = (s,0)$, and claim that
the function $F$ satisfies 
the Neumann boundary condition at $z$,
\[
\partial_t F (z) = 0 .
\]
Indeed, we compute at $z$ that
\begin{eqnarray*}  
\pp_t F 
\,=\,
df \left( \pp_t u \right) 
\,=\,
\langle \nabla f, \pp_t u \rangle
&=&
\langle X, \pp_t u \rangle \\
&=& \go \left( X, J \pp_t u \right)
\,=\,
- \go \left( X, \pp_s u\right) 
\,=\, 0 ,
\end{eqnarray*}
where in the last step we have used that $X \in TL_0$ by
(H4) and $\pp_s u \in TL_0$ by \eqref{floer}.
Let now $\tau$ be the reflection $(s,t) \mapsto (s,-t)$, 
set $\widehat{\Omega} = \Omega \cup \tau \left( \Omega \right)$,
and let $\widehat{F}$ be the extension of $F$ to
$\widehat{\Omega}$ satisfying $\widehat{F}(s, -t) = F(s,t)$.
Since $\pp_t F = 0$ along $\{ t=0 \}$, the continuous function
$\widehat{F}$ is weakly subharmonic, and hence cannot have a
strict maximum on $\widehat{\Omega}$.
Repeating this argument for $z = (s,1) \in \Omega$, we see that
the same holds for $F$ on $\Omega$, and so either $u(S) \cap V
= \emptyset$, or $F$ is locally constant. In the latter case,
$\Omega = S$, so that
\[
\lim_{s \ra \infty} u (s,t) = c_+ \in L_0 \cap V ,
\]
which is impossible in view of \eqref{e:empty}.
\proofend

We endow $\cj$ with the $C^\infty$-topology.
Recall that a subset of $\cj$ is {\it generic}\, if it is contained in
a countable intersection of open and dense subsets.
For $\bJ \in \cj$ let $\cm (\bJ)$ be the space of solutions of
\eqref{floer}.
The following proposition is proved in \cite{FHS, Oh2}.

\begin{proposition}
There exists a generic subset $\cj_{\reg}$ of $\cj$ 
such that for each $\bJ \in \cj_{\reg}$
the moduli space $\cm (\bJ)$ is a smooth finite 
dimensional manifold. 
\end{proposition}
Under hypotheses (H1)--(H4), the ungraded Floer homology $HF
\left( M,L_0,L_1 \right)$
can be defined. 
In order to prove Theorem~\ref{t:twist} we must compute the rank of this
homology, and to this end it will be crucial 
to endow it with a $\ZZ$-grading. We therefore
impose a final hypothesis.
For $\bJ \in \cj_{\reg}$ 
denote the submanifold of those $u \in \cm (\bJ)$ with
$\lim_{s \to \pm \infty}u(s,t)=c_\pm$ by 
$\cm (c_-,c_+;\bJ)$, and for $u \in \cm \left( c_-,c_+; \bJ \right)$ denote
by $I(u)$ the local dimension of $\cm \left( c_-,c_+; \bJ \right)$ at $u$.
\begin{itemize}
\item[\bf{(H5)}\,] 
$I(u)$ only depends on $c_-$ and $c_+$. 
\end{itemize}
Using (H5) and gluing one sees that there exists an index function
\[
\ind \colon L_0 \cap L_1 \to \ZZ
\]
such that $I(u) = \ind c_- - \ind c_+$,
so that
\[
\dim \cm(c_-,c_+;\bJ) \,=\, \ind c_- - \ind c_+ .
\]
For $k \in \ZZ$ let $CF_k (M,L_0,L_1)$ be the $\ZZ_2$-vector space
generated by the points $c \in L_0 \cap L_1$ with $\ind c = k$.
In view of (H1), the rank of $CF_k (M,L_0,L_1)$ is finite.
In order to define a chain map on $CF_* (M,L_0,L_1)$ we need the
following 
\begin{lemma}  \label{l:energy}
For $u \in \cm (c_-,c_+; \bJ)$ the energy $E(u)$ only depends on $c_-$
and $c_+$.
\end{lemma}
\proof
We have
$E(u) = \int_u d \gl = \int_{\pp u} \gl 
        = \int_{u (\RR,0)} \gl - \int_{u (\RR,1)} \gl$
for any primitive $\gl$ of $\go$.
Since $d \gl |_{L_j} =0$ and $H^1 \left( L_j; \RR \right) =0$, we find
smooth functions $f_i$ on $L_j$ such that $\gl |_{L_j} = d f_j$ for
$j=0,1$. Therefore, $E(u) = f_0 (c_+) - f_0 (c_-) - f_1 (c_+) +
f_1(c_-)$.
\proofend

The group $\RR$ acts on $\cm (c_-,c_+;\bJ)$ by time-shift.
In view of
Lemma~\ref{l:stayaway} the elements of $\cm (c_-,c_+;\bJ)$
uniformly stay away from the boundary $\pp M$,
and by Lemma~\ref{l:energy} and (H1), their energy is uniformly bounded.
Moreover, $[\go] |_{\pi_2(M)} =0$ and $[\go] |_{\pi_2(M,L_j)} =0$ since
$\go$ is exact and by (H3),
so that when taking limits in $\cm (c_-,c_+;\bJ)$ there is
no bubbling off of $\bJ$-holomorphic spheres or discs.
The Floer-Gromov compactness theorem
thus implies that the quotient $\cm (c_-,c_+;\bJ)/\RR$ 
is compact.
In particular, if $\ind c_- - \ind c_+ = 1$, then 
$\cm (c_-,c_+;\bJ)/\RR$ is a finite set, and we then set
\[
n(c_-,c_+;\bJ) \,=\, 
\# \left\{ \cm (c_-,c_+;\bJ)/\RR \right\} \mod 2 .
\]
For $k \in \ZZ$ define the Floer
boundary operator $\partial_k (\bJ) \colon CF_k \to CF_{k-1}$ as the
linear extension of
\[
\partial_k (\bJ) c \,=\, 
\sum_{\substack{c' \in L_0 \cap L_1\\i(c')=k-1}}
n(c',c)\,c' .
\]
Using the compactness of the $0$- and $1$-dimensional 
parts of $\cm (\bJ)/\RR$ one shows by gluing that
$\pp_{k-1} (\bJ) \circ \pp_k (\bJ)  = 0$ for each $k$, see \cite{F3,Sch1}.
The complex $\left( CF_* (M,L_0,L_1;\bJ), \pp_* (\bJ)\right)$ is called
the Floer chain complex. A continuation argument together with
Lemma~\ref{l:stayaway} shows that its homology
\[
HF_k(M,L_0,L_1;\bJ) \,=\,  \frac{\Ker \partial_k(\bJ)}{\Im \partial_{k+1}(\bJ)}  
\]
is a graded $\ZZ_2$-vector space which does not depend on 
$\bJ \in \cj_{\reg}$, see again \cite{F3,Sch1},
and so we can define the Lagrangian Floer homology of the
triple $(M,L_0,L_1)$ by
\[
HF_*(M,L_0,L_1) \,=\, HF_*(M,L_0,L_1;\bJ)
\]
for any $\bJ \in \cj_{\reg}$.
We denote by $\Hamc(M) \subset \Sympcc (M)$ the group of 
Hamiltonian diffeomorphisms generated by
time-dependent Hamiltonian functions whose support is contained in 
$\Int M$. The usual continuation argument also implies
\begin{proposition}  \label{p:Haminvariance}
For any $\gf \in \Hamc (M)$ we have 
$HF_* (M, \gf (L_0),L_1) = HF_* (M, L_0,L_1)$
as graded $\ZZ_2$-vector spaces.
\end{proposition}

\subsection{Computation of $HF_* \left( \gt^m (D_x), D_y \right)$}  
                                           \label{HF}
\begin{figure}[h] 
 \begin{center}
  \psfrag{c0+}{$c_0^+$}
  \psfrag{c0-}{$c_0^-$}
  \psfrag{c1+}{$c_1^+$}
  \psfrag{c1-}{$c_1^-$}
  \psfrag{g+}{$\dot{\gg}^+(\gd)$}
  \psfrag{g-}{$\dot{\gg}^-(\gd)$}
  \psfrag{Dy}{$D_y$}
  \leavevmode\epsfbox{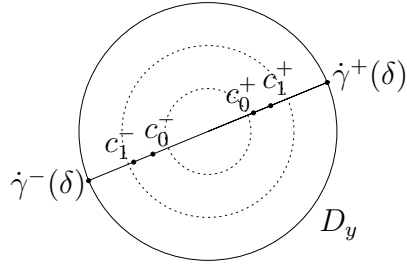}
 \end{center}
 \caption{The points $c_i^{\pm}  \in L_0 \cap L_1$ for $m=2$.}
 \label{figure3}
\end{figure}

\ni
We consider a $d$-dimen\-sio\-nal $P_k$-manifold $(B,g)$. 
Using the Riemannian metric $g$ we identify $T_1^*B$ with 
the unit ball bundle $T_1B$, and for $x \in B$ 
we set $D_x = T_xB \cap T_1B$.
According to Lemma~\ref{l:dicht} we find $x \in B$ such that 
$V_x = \left\{ v \in \pp D_x \mid l(v)=1 \right\}$ is a non-empty and
open subset of $\pp D_x$. 
We denote by $\rho$ the injectivity radius at $x$, and we define the
non-empty open subset $W$ of $B$ by
\[
W \,=\, \exp_x \left\{ v \in T_xB \mid 0 < \left| v \right| < \rho,\,
\frac{v}{|v|} \in V_x \right\} .
\]
Let $f \colon [0, \infty[\; \ra [0, \infty[$ be a smooth function as in
\eqref{e:f}.
More precisely, we choose $f$ such that
\[
f(r) = 0  \,\text{ if } r \in \left[ 0, \tfrac{1}{3} \right], \quad
f'(r) = 1 \,\text{ if } r \ge \tfrac{2}{3}, \quad
f''(r) >0 \,\text{ if } r \in \left] \tfrac{1}{3}, \tfrac{2}{3} \right[.
\]
Fix $m \in \ZZ \setminus \{ 0 \}$. 
For notational convenience we assume $m \ge 1$.
The symplectomorphism $\gt^m =
\gt_f^m \in \Symp^c \left( T^*B \right)$ is generated by $m f \left( |p|
\right)$.
Choose now $y \in W$.
By our choice of $f$ the two Lagrangian submanifolds
\[
L_0 = \gt^m (D_x) \quad \text{ and } \quad L_1 = D_y
\]
intersect transversely in exactly $2m$ points
and in particular meet hypothesis (H1); moreover, $\gt^m$ is the
identity on $U = T_1^* B \setminus T^*_{2/3}B$, so that
$L_0 \cap L_1 \cap U = \emptyset$ and (H2) is met.
Since $L_0$ and $L_1$ are simply connected, (H3) is also met, and
\begin{equation}  \label{def:X}
X \,=\, X(q,p) \,=\, \sum_{i=1}^d  p_i \frac{\partial}{\partial p_i}
\end{equation}
is a Liouville vector field defined on all of $T_1^*B$ which is
transverse to $\pp T_1^*B$, pointing outwards, and $X(x) \in T_x L_j$
for all $x \in L_j \cap U$, $j=0,1$, verifying (H4).
In order to verify (H5) we follow \cite{S1} and describe the natural
grading on $HF \left( T^*_1B, L_0, L_1 \right)$. 
Let $\gd$ be the distance of $y$ from $x$; then $0<\gd<\rho< 1/2$. 
For $i \in \NN_m = \left\{ 0,1, \dots, m-1 \right\}$ we set 
\[
\tau_i^+ = i+\gd 
\quad \text{ and } \quad 
\tau_i^- = i+1-\gd
\] 
and define $r_i^{\pm} \in \left] \tfrac{1}{3}, \tfrac{2}{3} \right[$
by
$m f' \left( r_i^{\pm} \right) = \tau_i^{\pm}$.
The $2m$ points in $L_0 \cap L_1$ are then given by
\[
c_i^\pm \,=\, \gt^m \left( r_i^{\pm} \dot{\gg}^{\pm} (0) \right) \,=\,
r_i^{\pm} \dot{\gg}^{\pm} (\gd) 
\]
where $\gg^+ \colon \RR \ra B$ is the geodesic with $\gg^+(0) = x$ and
$\gg(\gd) = y$
and $\gg^- (t) = \gg^+(-t)$ is the opposite geodesic, 
cf.\ Figure~\ref{figure3}.
%
%

\ni
Define the index function $\ind \colon L_0 \cap L_1 \ra \ZZ$ by 
\[
\ind c_i^\pm \,=\, \sum_{0<t< \tau_i^{\pm}} \ind \gg^{\pm} (t) , 
\quad \, i \in \NN_m .
\]
It is shown in \cite{S1} that for $\bJ \in \cj_{\reg}$ and $u \in \cm
\left( c_-, c_+; \bJ \right)$ the local dimension $I(u)$ of $\cm \left(
c_-, c_+; \bJ \right)$ at $u$ is $\ind c_- - \ind c_+$, so that (H5) is
also met.
We abbreviate 
\[
CF_* (B,m) = CF_* \left( T^*_1B, L_0, L_1 \right)
\quad \text{and} \quad
HF_* (B,m) = HF_* \left( T^*_1B, L_0, L_1 \right).
\]
Our next goal is the compute the Floer chain groups $CF_* (B,m)$.
\begin{proposition}  \label{p:jacobiper}
Let $(B,g)$ be a $d$-dimensional $P_k$-manifold and $i \in \NN_m$. Then
\[
\ind c_i^+ = i (k+d-1) 
\quad \text{ and } \quad
\ind c_i^- = i (k+d-1) +k . 
\]
\end{proposition} 

\proof
We start with a general
\begin{lemma}  \label{l:Jacobi}
Let $\gg \colon \RR \ra B$ be a geodesic of a $P$-manifold $(B,g)$, and
let $J \colon \RR \ra B$ be a Jacobi field along $\gg$
such that $J(0) = 0$.
If $J(t) = 0$, then $J(t+n) =0$ for all $n \in \ZZ$. 
\end{lemma}

\proof
We can assume that $J'(t) \neq 0$ since otherwise $J \equiv 0$.
Fix $n \in \ZZ$.
Since $(B,g)$ is a $P$-manifold, $\gg (t+n)$ is conjugate to $\gg (t)$
with multiplicity $d-1$.
Since this is the maximal possible multiplicity of a conjugate point,
and since $J(t) =0$ and $J'(t) \neq 0$, 
$J$ must be a Jacobi field conjugating $J(t)$ and $J(t+n)$, i.e.,
$J(t+n) =0$.
\proofend

By our choice of $y$, the point $y$ is not conjugate to $x$ along
$\gg^+ \colon [0,\gd] \ra B$, and so $\ind c_0^+ =0$.
Since $(B,g)$ is a $d$-dimensional $P_k$-manifold, Lemma~\ref{l:Jacobi}
now implies $\ind c_i^+ = i (k+d-1)$ for all $i \in \NN_m$.
The choice of $y$ and Lemma~\ref{l:Jacobi} imply that $\ind c_0^- = k$,
and now Lemma~\ref{l:Jacobi} implies $\ind c_i^- = i (k+d-1) +k$
for all $i \in \NN_m$.
\proofend

In view of Proposition~\ref{p:jacobiper} we find
\begin{corollary}  \label{c:grading}
Let $(B,g)$ be a $d$-dimensional $P_k$-manifold.
\begin{align*}
\intertext{If $k \ge 1$,} 
CF_i(B,m) &\,=\, \left\{\begin{array}{cl}
                   \ZZ_2 & i \in (k+d-1) \NN_m \cup 
                   \big( (k+d-1)\NN_m +k \big) \\
                   0 & \mathrm{otherwise};
                \end{array}\right.
\intertext{if $k=0$ and $d>1$,}
CF_i(B,m) &\,=\, \left\{\begin{array}{cl}
                   \mathbb{Z}_2^2 & i \in (d-1) \NN_m\\
                   0 & \mathrm{otherwise};
                \end{array}\right.
\intertext{if $k=0$ and $d=1$,}
CF_i(B,m) &\,=\, \left\{\begin{array}{cl}
                   \mathbb{Z}_2^{2m} & i=0\\
                   0 & \mathrm{otherwise}.
                \end{array}\right.
\end{align*}
\end{corollary}
\begin{theorem}  \label{t:floer} 
Assume that $(B,g)$ is a $P_k$-manifold as in Theorem~\ref{t:twist},
and if $k=1$ assume that $(B,g) = \left( \CP^n,g_{\can} \right)$.
Then the Floer boundary operator $\pp_* \colon CF_* (B,m) \ra CF_{*-1} (B,m)$
vanishes identically, and so $HF_*(B,m) = CF_*(B,m)$.
In particular, $\rank HF (B,m) = 2m$.
\end{theorem}

\proof
We first assume that $k \neq 1$ and $d \neq 2$.
Corollary~\ref{c:grading} then shows that for any $* \in \ZZ$, at least
one of the chain groups $CF_* (B,m)$ and $CF_{*-1} (B,m)$ 
is trivial, and so $\pp_* =0$.
It remains to prove the vanishing of $\pp_*$ for the spaces
$\RP^2$ and $\CP^n$, $n \ge 1$, endowed with their canonical $P$-metrics.
We shall do this by using a symmetry argument. 

\s
\ni
{\bf The case $\left( \CP^n, g_{\can} \right)$.}
Note that every
diffeomorphism $\gf$ of $\CP^n$ lifts to a symplectomorphism
$(\gf^{-1})^*$ of $T^*\CP^n$, and that if $\gf$ is an isometry,
then $(\gf^{-1})^*$ is a symplectomorphism of $T^*_1\CP^n$. 
Let $\RP^n$ be the real locus of $\CP^n$.
\begin{lemma}  \label{l:RPn}
We can assume without loss of generality that $x, y \in \RP^n$.
\end{lemma}

\proof
Choose a unitary matrix $U \in \mathbf{U} (n+1)$ such that $x' = U(x)
= [1:0: \dots : 0] \in \RP^n$ and $y' = U(y) \in \RP^n$. We again
identify $T^*_1\CP^n$ with $T_1\CP^n$ via the Riemannian metric $g_{\can}$.
Since $U$ is an isometry of $\left( \CP^n, g_{\can} \right)$, its lift $U_*$ to
$T_1\CP^n$ commutes with the geodesic flow on $T_1\CP^n$, and hence
$(U^{-1})^*$ commutes with $\gt^m$. 
Therefore,
\[
(U^{-1})^* L_0 \,=\, (U^{-1})^* \gt^m \left( D_x \right) \,=\, 
\gt^m (U^{-1})^* \left( D_x \right) \,=\,  \gt^m \left( D_{x'} \right)
\]
and $(U^{-1})^* L_1 = (U^{-1})^* D_y = D_{y'}$.
By the natural invariance of Lagrangian Floer homology we thus obtain
\begin{eqnarray*}
HF_* \left( T_1 \CP^n, L_0,L_1 \right) 
&=&
HF_* \left( T_1 \CP^n, (U^{-1})^*L_0,(U^{-1})^*L_1 \right)  \\
&=&
HF_* \left( T_1 \CP^n, \gt^m \left( D_{x'} \right),D_{y'} \right) ,
\end{eqnarray*}
as desired.
\proofend

Consider the involution
\begin{equation}  \label{inv:complex}
[ z_0 : z_1 : \dots : z_n ] \,\mapsto\,
[ \bar{z}_0 : \bar{z}_1 : \dots : \bar{z}_n ]
\end{equation}
of $\CP^n$.
Its fixed point set is $\RP^n$. Since complex conjugation
\eqref{inv:complex} is an isometry of
$\left( \CP^n, g_{\can} \right)$, it lifts to a
symplectic involution $\gs$ of $T^*_1 \CP^n$.
Since $x,y \in \RP^n$ and since complex conjugation is an isometry,
we see as in the proof of Lemma~\ref{l:RPn} that $\gs (L_j)
=L_j$, $j=0,1$,
and $\sigma$ acts trivially on $L_0 \cap L_1$. 
Assume that $\bJ \in \cj \left( T^*_1 \CP^n \right)$ is invariant under $\sigma$, 
i.e., $\sigma^* J_t = \gs_* J_t \gs_* = J_t$ for every $t \in [0,1]$. 
Then $\sigma$ induces an involution on the solutions of \eqref{floer} by
\[
u \mapsto \gs \circ u.
\]
If $u$ is invariant under $\sigma$, i.e., $u=\sigma \circ u$,
then $u$ is a solution of \eqref{floer} with 
$M$ replaced by the fixed point set $M^\gs = T^*_1\RP^n$ of $\gs$ 
and $L_j$ replaced by $L_j^\sigma=L_j \cap M^\sigma$ for $j =0,1$.
According to Proposition~\ref{p:index}, $\CP^n$ is a $P_1$-manifold and
$\RP^n$ is a $P_0$-manifold, and so we read off from
Corollary~\ref{c:grading} that if $\ind_M (c_-) - \ind_M (c_+) =1$, then 
$\ind_{M^\sigma}(c_-) - \ind_{M^\sigma}(c^+) =0$. 
One thus expects that
for generic $\sigma$-invariant $\bJ \in \cj$
there are no solutions of \eqref{floer} which are invariant under $\sigma$.
In particular, solutions of \eqref{floer} appear in pairs, and so
$\pp_*=0$. To make this argument precise, we
need to show that there exist
$\sigma$-invariant $\bJ \in \cj$
which are ``regular'' for every non-invariant solution of \eqref{floer} and
whose restriction to $M^\sigma$ is also ``regular''. This will be done
in the next paragraph.

\subsubsection[A transversality theorem]{A transversality theorem}
                              \label{transverse}

\ni
We consider, more generally, an exact compact symplectic manifold
$(M,\go)$ with boundary $\pp M$ containing two compact
Lagrangian submanifolds $L_0$ and $L_1$ as in \ref{ss:LagFloer}:
(H1), (H2), (H3) hold and there is a Liouville vector field
$X$ on a neighbourhood $U$ of $\pp M$ such that (H4) holds.
We in addition assume that $\gs$ is a symplectic involution of
$(M,\go)$ such that
\begin{equation}  \label{e:I123}
\gs (L_j) = L_j,\,\, j=0,1,
\qquad
\gs |_{L_0 \cap L_1} = \id,
\qquad
\gs_* X = X .
\end{equation}
We have already verified the first two properties for $M = T_1^*
\CP^n$ and the lift $\gs$ of \eqref{inv:complex}, and we notice
that $\gs_* X=X$ for the Liouville vector field \eqref{def:X}.
The fixed point set $M^\gs = \mathrm{Fix}(\gs)$ is a symplectic
submanifold of $(M, \go)$. 
Set $\go^\gs = \go |_{M^\gs}$. Since $\gs_* X=X$ the vector
field $X^\gs = X |_{U \cap M^\gs}$ is a Liouville vector field
near $\pp M^\gs$.
As in \ref{ss:LagFloer} we denote by $\cj = \cj (M)$ the space of
smooth families $\bJ = \{ J_t\}$, $t \in [0,1]$, 
of smooth $\go$-compatible almost complex structures on $M$
which on $V$ do not depend on $t$
and meet (J1), (J2), (J3).
The space $\cj \left( M^\gs \right)$ is defined analogously by
imposing (J1), (J2), (J3) for $X^\gs$ on $M^\gs \cap V$.
The subspace of those $\bJ$ in $\cj (M)$ which are
$\gs$-invariant is denoted $\cj^\gs (M)$.
There is a natural restriction map
\[
\varrho \colon \cj^\gs (M) \to \cj \left(M^\gs \right), \quad\,
\bJ \mapsto \bJ|_{TM^\gs}.
\]
\begin{lemma}  \label{l:rho}
The restriction map $\varrho$ is open.
\end{lemma}

\proof 
Recall that $\gf_r$, $r \le 0$, denotes the semiflow of $X$, and
that $\xi$ and $R$ are the contact structure and the Reeb vector
field on $\pp M$ defined by \eqref{def:xi} and
\eqref{def:Reeb}. Since $\gs$ is symplectic, $\gs_* X =X$ and
$\gs (\pp M) = \pp M$ we have 
\[
\gs_* \xi = \xi 
\qquad \text{ and } \qquad
\gs_* R =R .
\]
The contact structure $\xi^\gs$ on $\pp M^\gs$ associated with
$X^\gs$ is 
$\xi \cap T \pp M^\gs$, and the Reeb vector field $R^\gs$ is $R
|_{\pp M^\gs}$.
We shall prove Lemma~\ref{l:rho} by first showing that $\rho$ is
onto. From the proof it will then easily follow that $\rho$ is open.

\s
\ni
{\bf Step 1. $\rho$ is onto:}
Fix $\bJ^\gs \in \cj \left(M^\gs \right)$. 
We set $g_t^\gs = \go \circ \left( \id \times J_t^\gs \right)$.
Choose a smooth family $\bg = \{ g_t \}$, $t \in [0,1]$, of
Riemannian metrics on $TM$ which on $V$ does not depend on $t$
and satisfies
\begin{itemize}
\item[(g1)] 
$g \left( X(x), v \right) = 0, 
\qquad 
x \in \pp M, \,\, v \in T_x \pp M,$ \vspace{0.4em}
\item[(g2)] 
$g \left( X(x), X(x) \right) =1,
\qquad 
x \in \pp M ,$ \vspace{0.4em}
\item[(g3)] 
$\gf^*_r g (x) = e^r g(x), 
\qquad 
x \in \pp M, \,\, r \in [-\eps', 0] ,$ \vspace{0.4em}
\item[(g4)] 
$g \left( R(x), R(x) \right) = 1, 
\qquad
x \in \pp M ,$ \vspace{0.4em}
\item[(g5)]
$g \left( R(x), v \right) = 0,
\qquad
x \in  \pp M , \,\, v \in \xi ,$
\end{itemize}
and in addition satisfies for each $t$
\begin{itemize}
 \item[(g6)] 
  if $x \in M^\gs$, 
  then $g_t(x)|_{T_x M^\gs} = g^\gs_t (x)$, \vspace{0.4em}
 \item[(g7)]
  if $x \in M^\gs$, 
  then the Riemannian and the symplectic
  orthogonal complement of $T_x M^\gs$ in $T_x M$ 
  coincide, i.e.,
  if for $\eta \in T_x M$ it holds that 
  $\go (\eta,\zeta)=0$ for 
  every $\zeta \in T_x M^\gs$, then also 
  $g(\eta,\zeta)=0$ for
  every $\zeta \in T_x M^\gs$, \vspace{0.4em}
 \item[(g8)] $\gs^* g_t = g_t$.
\end{itemize}
In order to see that such a family $\bg$ exists, first notice that
in view of (J1), (J2), (J3), Lemma~\ref{l:equiv} and
\eqref{def:Reeb},
the metric $g^\gs$ satisfies (g1)--(g5) for $X^\gs$, $R^\gs$,
$x\in \pp M^\gs$ and $v \in T_x\pp M^\gs$ or $v \in \xi^\gs$.
We thus find a family $\bg_0$ satisfying (g1)--(g7).
Then $\gs^* \bg_0$ also satisfies (g1)--(g7) as one readily
verifies; we only mention that (g3) follows from $\gs^* \circ \gf_r =
\gf_r \circ \gs$ which is a consequence of $\gs_* X=X$. Now set
\[
\bg \,=\, \tfrac{1}{2} \left( \bg_0 + \gs^* \bg_0 \right) .
\]  
Let $\mathfrak{Met}$ be the space of smooth Riemannian metrics on
$M$ and let $\cj (\go)$ be the space of smooth $\go$-compatible
almost complex structures on $M$.
For $J \in \cj (\go)$ we write $g_J = \go \circ \left( \id
\times J \right) \in \mathfrak{Met}$.
It is shown in \cite[Proposition 2.50\:(ii)]{MS} 
that there exists a smooth map 
\[
r \colon \mathfrak{Met} \to \cj (\go),
\quad \,
g \mapsto r(g) =: J_g
\]
such that
\begin{equation}  \label{prop:R}
r \left( g_J \right) = J
\qquad \text{ and } \qquad 
r \left( \gf^* g \right) = \gf^* r(g)
\end{equation}
for every symplectomorphism $\gf$ of $M$. 
We define $\bJ = \{ J_t \}$ by
\[
J_t \,=\, r \left( g_t \right) .
\]
The second property in \eqref{prop:R} and (g8) show that $\bJ$
is $\gs$-invariant.
In order to prove that $\bJ \in \cj^\gs (M)$ we also need to
show that each $J_t$ meets (J1), (J2), (J3) and $J_t
|_{M^\gs} = J_t^\gs$.
To this end we must recall the construction of $r$ from \cite{MS}. 
Fix $g \in \mathfrak{Met}$ and $x \in M$. The automorphism $A$
of $T_xM$ defined by $\go_x (v,w)= g_x(Av,w)$ is
$g$-skew-adjoint. Denoting by $A^*$ its $g$-adjoint, we find
that $P = A^*A = -A^2$ is $g$-positive definite.
Let $Q$ be the automorphism of $T_xM$ which is $g$-self-adjoint,
$g$-positive definite, and satisfies $Q^2=P$, and then set 
\[
J_x (\go,g) \,=\, Q^{-1}A . 
\]
It is clear that $J_x(\go,g)$ depends smoothly on $x$.
The map $r$ is defined by $r(g)(x) = J_x(\go,g)$.
One readily verifies that $r(g)$ is $\go$-compatible, see
\cite[p.\ 14]{HZ}, and meets \eqref{prop:R}.
From the construction we moreover read off that
\begin{itemize}
 \item[(r1)] 
  $J_x ( c_1 \go, c_2 g ) = J_x(\go,g)$ for all $c_1,c_2>0$, \vspace{0.4em}
 \item[(r2)]
  if $T_xM =V \oplus W$ in such a way that $W$ is both
$\go$-orthogonal and $g$-orthogonal to $V$, i.e., $W=V^\go =
V^\perp$, then $A$, $P$ and $Q$ leave both $V$ and $W$
invariant, so that $J_x(\go,g)$ leaves $V$ invariant and 
\[
J_x(\go,g) |_V \,=\, J_x \left( \go |_V, g |_V \right) .
\]
\end{itemize}
We are now in position to verify (J1), (J2), (J3) for $J_t =
r(g_t) = J_{g_t}$.
In view of \eqref{def:xi} and \eqref{def:Reeb} and (g1) and (g5)
the plane field $\langle X,R \rangle$ on $\pp M$ generated by
$X$ and $R$ is both $\go$-orthogonal and $g$-orthogonal to
$\xi$,
\[
\langle X,R \rangle \,=\, \xi^\go \,=\, \xi^\perp ,
\]
and so (r2) implies 
\begin{equation}  \label{e:Jg}
J_g | \langle X,R \rangle \,=\, J_{g | \langle X,R \rangle} .
\end{equation}
Define the complex structure $J_0$ on $\langle X,R \rangle$ by
$J_0 X=R$. Using (g1), (g2), (g4) and \eqref{def:Reeb} we find 
$g | \langle X,R \rangle = g_{J_0}$, and so the first property
in \eqref{prop:R} implies $J_{g | \langle X,R \rangle} =
J_0$. Together with \eqref{e:Jg} we find 
\begin{equation}  \label{e:g0}
J_g | \langle X,R \rangle \,=\, J_0 . 
\end{equation}
The $J_g$-invariance of $\go$, \eqref{e:g0} and \eqref{def:Reeb} yield
(J1) and (J2).
The identity (J3) follows from $\gf_r^* \go = e^r \go$, (g3) and
(r1). 
Finally, $J_t |_{M^\gs} = J_t^\gs$ follows from (g6), (g7), (r2)
and the first property in \eqref{prop:R}.

\s
\ni
{\bf Step 2. $\rho$ is open:}
Let $U$ be an open subset of $\cj^\gs (M)$.
We must show that given $\bJ^\gs \in \rho (U)$, every ($C^\infty$)-close
enough $\widetilde{\bJ}^\gs \in \cj \left( M^\gs \right)$ belongs to $\rho
(U)$.
Fix $\bJ \in U$ with $\rho (\bJ) =
\bJ^\gs$, and set $\bg = g_{\bJ}$.
Since $\bJ \in \cj^\gs (M)$, the family $\bg$ satisfies
(g1)--(g8).
If $\widetilde{\bJ}^\gs \in \cj \left( M^\gs \right)$ is close
to $\bJ^\gs$, then 
$\widetilde{\bg}^\gs = g_{\widetilde{\bJ}^\gs}$  
is close to $g_{\bJ^\gs}$, and so we can choose
a smooth family $\widetilde{\bg}_0$ close to $\bg$ which
satisfies (g1)--(g7).
Then
\[
\widetilde{\bg} \,=\, 
\tfrac 12 \left( \widetilde{\bg}_0 + \gs^*\widetilde{\bg}_0 \right)
\]
satisfies (g1)--(g8), and since $\widetilde{\bg}_0$ was close to
$\bg$ and since $\gs^* \bg = \bg$, the family $\widetilde{\bg}$
is also close to $\bg$.
Set $\widetilde{\bJ} = r \left( \widetilde{\bg} \right)$.
Then $\rho \big( \widetilde{\bJ} \big) =
\widetilde{\bJ}^\gs$, and since $r \colon \mathfrak{Met} \to
\cj (\go)$ is smooth and $\widetilde{\bg}$ is close to $\bg$, we
see that $\widetilde{\bJ} = r \left( \widetilde{\bg} \right)$ is close
to $r (\bg) = r \left( g_{\bJ} \right) = \bJ$.
In particular, if $\widetilde{\bJ}^\gs$ was close enough to
$\bJ^\gs$, then $\widetilde{\bJ} \in U$. 
The proof of Lemma~\ref{l:rho} is complete. 
\proofend

For the remainder of the proof of Theorem~\ref{t:floer} for 
$\left( \CP^n, g_{\can} \right)$ we assume that the reader is
familiar with the standard transversality arguments in Floer
theory as presented in Section~5 of \cite{FHS} or Sections~3.1 and 3.2 of
\cite{MS3}, 
and we shall focus on those aspects of the argument particular to
our situation. 
Fix $c_-,c_+ \in L_0 \cap L_1$.
We interpret solutions of \eqref{floer} with $\lim_{s \ra
\pm \infty} u(s,t) = c_\pm$ as the zero set of a smooth section
from a Banach manifold $\cb$ to a Banach bundle $\ce$ over $\cb$.
We fix $p>2$. 
According to Lemma~D.1 in \cite{RS} there exists a smooth family
of Riemannian metrics $\{ g_t \}$, $t \in [0,1]$, on $M$ such
that $L_j$ is totally geodesic with respect to $g_j$, $j=0,1$.
Let $\cb = \cb^{1,p} \left(c_-,c_+ \right)$ be the space of
continuous maps $u$ from the strip $S = \RR \times [0,1]$ 
to the interior $\Int M$ of $M$
which satisfy 
$\lim_{s \ra \pm \infty} u(s,t) = c_\pm$ uniformly in $t$,
are locally of class $W^{1,p}$, and satisfy the conditions
\begin{itemize}
\item[(B1)] $u(s,j) \in L_j$ for $j =0,1$, \vspace{0.2em}
\item[(B2)]  
  there exists $T>0$,
  $\xi_- \in W^{1,p}((-\infty,-T] \times [0,1],T_{c_-}M)$, and
  $\xi_+ \in W^{1,p}([T,\infty) \times [0,1], T_{c_+}M)$ with
  $\xi_\pm(s,j) \in T_{c_\pm}L_j$ such that
\begin{eqnarray*}
u(s,t)&=& 
   \left\{ \begin{array}{ll} 
      \exp_{c_-}(\xi_-(s,t)), & s \leq -T, \\
      \exp_{c_+}(\xi_+(s,t)), & s \geq T.
           \end{array}
   \right.
\end{eqnarray*}
\end{itemize}
Here, $\exp_{c_\pm} \left( \xi_\pm(s,t) \right)$ denotes the
image of $\xi_\pm(s,t)$ under the exponential map with respect
to $g_t$ at $c_\pm$.
The space $\cb$ is an infinite dimensional Banach manifold 
whose tangent space at $u$ is 
\[
T_u\cb \,=\, \left\{ \xi \in W^{1,p} \left( S,u^*TM \right) \mid
 \xi(s,j) \in T_{u(s,j)}L_j,\, j=0,1 \right\} .
\]
Let $\ce$ be the Banach bundle over $\cb$ whose fibre over $u \in \cb$ is 
\[
\ce_u \,=\, L^p \left( S,u^*TM \right).
\]
For $\bJ \in \cj (M)$ define the section $\cf_\bJ \colon \cb \to
\ce$ by
\[
\cf_{\bJ}(u) \,=\, \pp_s u + J_t(u) \pp_t(u) 
\]
and set $\cm_\bJ = \cf^{-1}_\bJ(0)$.
The set $\cm_\bJ$ agrees with the set of those $u \in \cm (\bJ)$
with $\lim_{s \ra \pm \infty} u(s,t) = c_\pm$. 
Indeed, Lemma~\ref{l:stayaway} and Proposition~1.21 in \cite{FHS}
show that the latter set belongs to $\cm_\bJ$.
Conversely, in view of $p>2$, 
elliptic regularity and (B2) imply that $u \in \cm_\bJ$ is
smooth and satisfies 
$\lim_{s \ra \pm \infty} \pp_s u(s,t) = 0$ uniformly in $t$, 
so that $E(u)<\infty$ by Proposition~1.21 in \cite{FHS}. 
If $u \in \cm_\bJ$, then the vertical differential
of $\cf_\bJ$, 
\[
D_{u,\bJ} \equiv D \cf_\bJ(u) \colon T_u \cb \to \ce_u, \quad
 \xi \mapsto \nabla_s \xi +J_t(u) \nabla_t \xi + \nabla_\xi J_t(u)
\pp_t u ,
\]
is a Fredholm operator, cf.\ \cite[Theorem~2.2]{Sa}.
Here, $\nabla$ denotes the Levi--Civita connection with respect
to the $t$-dependent metric $g_{J_t}$.
We further consider the Banach submanifold
\[
\cb^\sigma \,=\, \left\{ u \in \cb \mid u = \sigma \circ u \right\},
\]
of those $u$ in $\cb$ whose image lies
in $M^\sigma$. We denote by $\ce^\sigma$ the Banach bundle over
$\cb^\sigma$ whose fibre over $u \in \cb^\sigma$ is
\[
\ce^\gs_u \,=\, L^p \left( S, u^*TM^\gs \right) .
\]
Note that $\ce^\gs$ is a subbbundle of the restriction of
$\ce$ to $\cb^\gs$.
For $\bJ \in \cj^\gs (M)$ we abbreviate
\[
\cm^\gs_\bJ \,\equiv\, \cf_\bJ^{-1}(0) \cap \cb^\gs \,=\,
\cm_\bJ \cap \cb^\gs 
\]
and for $u \in \cm^\gs_\bJ$ we set
\[
D^\gs_{u,\bJ} \equiv D_{u,\bJ} |_{T_u \cb^\gs} \colon 
T_u\cb^\gs \to \ce^\gs_u.
\]
\begin{definition}
{\rm
We say that $\bJ \in \cj^\gs (M)$ is {\it regular} if for every
$u \in \cm_\bJ \setminus \cm^\sigma_\bJ$ the operator
$D_{u,\bJ}$ is onto and if for every $u \in \cm^\sigma_\bJ$ the
operator $D^\sigma_{u,\bJ}$ is onto. 
}
\end{definition}
\begin{proposition}  \label{p:transverse}
The set $\left( \cj^\gs (M) \right)_{\reg}$ of regular almost complex
structures is generic in $\cj^\gs (M)$. 
\end{proposition}

\proof
It is proved in \cite[Proposition~5.13]{KS} that the subset
$\cR_1$ of those $\bJ \in \cj^\gs (M)$ for which $D_{u,\bJ}$ is onto for
every $u \in \cm_\bJ \setminus \cm^\sigma_\bJ$ is generic
in $\cj^\gs (M)$. 
Moreover, it is proved in \cite[Section~5]{FHS} that the subset
$\cR_2^\gs$ of those $\bJ^\gs \in \cj \left( M^\gs \right)$ for
which $D_{u,\bJ^\gs}$ is onto for every $u \in \cm_{\bJ}^\gs$ is
generic in $\cj \left( M^\gs \right)$. Notice that for $\bJ \in
\cj^\gs (M)$ we have $\cm_{\bJ}^\gs = \cm_{\rho ( \bJ )}$ and
$D_{u,\bJ}^\gs = D_{u,\rho (\bJ)}$ for $u \in \cm_\bJ^\gs =
\cm_{\rho (\bJ)}$. It follows that for $\bJ \in \cR_2 \equiv
\rho^{-1} \left( \cR_2^\gs \right)$ the operator
$D_{u,\bJ}^\gs$ is onto for every $u \in \cm_{\bJ^\gs}$. Since
the preimage of a generic set under a continuous open map is
generic, $\cR_2$ is generic in $\cj^\gs (M)$. Therefore, the set
of regular $\bJ \in \cj^\gs(M)$ contains the generic set $\cR_1
\cap \cR_2$, and the proof of Proposition~\ref{p:transverse} is
complete.
\proofend

In order to complete the proof of Theorem~\ref{t:floer} for 
$\left( \CP^n, g_{\can} \right)$, set again $M = T_1^*
\CP^n$. In view of Proposition~\ref{p:transverse} we find a $\bJ
\in \cj^\gs (M)$ which is regular for all $c_- , c_+ \in L_0 \cap
L_1$.
Fix $c_- , c_+$ with $\ind_M (c_-) - \ind_M (c_+)=1$. Since 
$\ind_{M^\gs} (c_-) - \ind_{M^\gs} (c_+)=0$, the Fredholm index
of $D_{u,\bJ}^\gs$ for $u \in \cm_\bJ^\gs$ vanishes, so that the
manifold of solutions of \eqref{floer} contained in $M^\gs$ is
$0$-dimensional and hence empty. 
Moreover, $D_{u,\bJ}$ is onto for every $u \in \cm_\bJ \setminus
\cm_\bJ^\gs$, and so $D_{u,\bJ}$ is onto for every $u \in
\cm_\bJ$. We can thus compute the Floer homology 
$HF_* \left( M, L_0,L_1 \right)$ by using $\bJ$. Since
$\cm_\bJ^\gs$ is empty, $\pp_* =0$, and the proof of Theorem~\ref{t:floer} for 
$\left( \CP^n, g_{\can} \right)$ is complete.

\s
\ni
{\bf The case $\left( \RP^2, g_{\can} \right)$.}
We consider the involution 
\begin{equation}  \label{inv:RP2}
[ q_0 : q_1 : q_2 ] \,\mapsto\, [ -q_0 : q_1 : q_2 ]
\end{equation}
of $\RP^2$. Its fixed point set is $[1:0:0] \cup \{
[0:q_1:q_2]\} \equiv p_N \cup \RP^1 = p_N \cup S^1$.
Since every isometry of $\left( S^2, g_{\can} \right)$ descends
to an isometry of $\left( \RP^2, g_{\can} \right)$, we find an
isometry of $\RP^2$ mapping $x$ and $y$ to $S^1$;
by the argument in the proof of Lemma~\ref{l:RPn} we can thus
assume that $x,y \in S^1$.
Let $\gs$ be the symplectic involution of $T_1^* \RP^2$ obtained
from lifting \eqref{inv:RP2}.
It satisfies \eqref{e:I123}.
The fixed point set of $\gs$ is $p_N \cup T_1^*S^1$.
Since $p_N$ is disjoint from $L_0 \cap L_1$, there is no
solution of \eqref{floer} lying in $p_N$. Set $M^\gs = T^*_1S_1$.
According to Proposition~\ref{p:index}, 
both $\RP^2$ and $S^1$ are $P_0$-manifolds, and so 
Corollary~\ref{c:grading} implies that if 
$\ind_M (c_-) - \ind_M (c_+) =1$, then 
$\ind_{M^\sigma}(c_-) - \ind_{M^\sigma}(c^+) =0$. 
The vanishing of $\pp_*$ now follows as in \ref{transverse}.
The proof of Theorem~\ref{t:floer} is finally complete.
\proofend

\subsection{End of the proof}  \label{end3}

\ni
For $(B,g) = S^1$, Theorem~\ref{t:twist} follows from the topological
argument given in Section~\ref{idea}, see also Corollary~\ref{c:top} below.
For the remainder of this subsection we therefore assume that $(B,g)$ is a
$P$-manifold of dimension $d \ge 2$.
We abbreviate $M=T^*B$ and $M_r = T_r^*B$.

\begin{lemma}  \label{l:Ham}
$\Hamc (M)= \Sympc0 (M)$.
\end{lemma}
\proof
Since $T^*B$ is orientable, Poincar\'e duality yields
\[
H_c^1 \left( M; \RR \right) \,\cong\,
H_{2d-1} \left( M; \RR \right) \,\cong\,
H_{2d-1} \left( B; \RR \right) \,=\, 0 .
\]
The lemma now follows in view of the exact sequence~\eqref{s:exact}.
\proofend

Let now $(B,g)$ be a $P$-manifold as in Theorem~\ref{t:floer}.
Let $\gt = \gt_f$ be the twist considered in Subsection~\ref{HF}, 
and let $\gf \in \Sympcc (M)$ be such that
$[\gf] = \left[ \gt^m \right] \in \pi_0 \left( \Sympcc (M) \right)$
for some $m \in \ZZ \setminus \{0\}$.
We assume without loss of generality that $m \ge 1$.
By Lemma~\ref{l:Ham} we find $r>0$ such that 
$\gt^m \gf^{-1} \in \Ham^c \left(M_r \right)$.
Then $\gt^{mn} \gf^{-n} \in \Ham^c \left(M_r \right)$ for all $n \ge 1$.
We assume without loss of generality that $r=1$.
Let $W$ be the non-empty open subset of $B$ constructed in \ref{HF}, and
fix $y \in W$.
We first assume that $\gf^n (D_x)$ intersects $D_y$ transversally.
Then $HF \left( M_1, \gf^n (D_x), D_y \right)$ is defined, and
in view of Proposition~\ref{p:Haminvariance} and Theorem~\ref{t:floer} 
we find that  
\begin{eqnarray*}
\rank CF \big( M_1, \gf^n \left( D_x \right), D_y \big) &\ge& 
\rank HF \big( M_1, \gf^{n} \left( D_x \right), D_y \big) \\
     &=& \rank HF  \big( M_1, \gt^{mn} \left( D_x \right), D_y \big)  \\
     &=& 2mn .
\end{eqnarray*}
It follows that the $d$-dimensional submanifold $\gf^n(D_x)$ of $M_1$
intersects $D_y$ at least $2mn$ times. Since this holds true for
every $y \in W$ and since $\pi \colon \left( M_1, g^* \right) \ra
(B,g)$ is a Riemannian submersion, we conclude that
\begin{equation*}  
\mu_{g^*} \big( \gf^n (D_x) \big) \,\ge\, 2mn \mu_g(W) .
\end{equation*}
If $\gf^n(D_x)$ and $D_y$ are not transverse, we choose a sequence
$\gf_i \in \Symp^c (M_1)$ such that $\gf_i^n (D_x)$ and $D_y$ are
transverse for all $i$, and $\gf_i \ra \gf$ in the $C^\infty$-topology.
For $i$ large enough, 
$\left[ \gf_i \right] = [ \gf ] \in \pi_0 \left(\Symp^c(M_1) \right)$, and
\begin{equation}  \label{est:fU}
\mu_{g^*} \big( \gf^n (D_x) \big) \,=\,
\lim_{i \ra \infty} \mu_{g^*} \big( \gf_i^n (D_x) \big) \,\ge\, 
2mn \mu_g(W) .
\end{equation}
Choose a smooth embedding $\gs \colon Q^d \ra T_x^*B$ such that $D_x
\subset \gs (Q^d)$. Then 
$\mu_{g^*} \left( \gf^n (\gs) \right) \ge \left( 2 m \mu_g(W) \right) n$, 
and so $s_d(\gf) \ge l (\gf) \ge 1$, as claimed.

\s
Assume next that $(B,g) = \left( \CP^{2n-1} / \ZZ_2, g_{\can} \right)$
and that $\left[ \gf \right] = \left[ \gt^m \right] \in \pi_0
\left(\Symp^c(M) \right)$.
Set $\big( \widetilde{B}, \tilde{g} \big) = \left( \CP^{2n-1},
g_{\can} \right)$ and $\widetilde{M} = T^* \widetilde{B}$.
In view of our normalization of $g$ and $\tilde{g}$, the twist
$\tilde{\gt}$ on $\widetilde{M}$ is a lift of $\gt$, and so
$\tilde{\gt}^m$ is a lift of $\gt^m$.
Lifting a symplectic isotopy between $\gt^m$ and $\gf$ to $\widetilde{M}$,
we obtain a symplectic isotopy between $\tilde{\gt}^m$ and a lift
$\tilde{\gf}$ of $\gf$. Since the projection $\widetilde{M} \ra M$ is a
local isometry, we thus obtain from \eqref{est:fU} that
\[
\mu_{g^*} \big( \gf^n (D_x) \big) \,\ge\,
\tfrac 12 \mu_{\tilde{g}^*} \big( \tilde{\gf}^n (D_{\tilde{x}}) \big)
\,\ge\, mn \mu_{\tilde{g}}(W) 
\]
for any $x \in B$ and a lift $\tilde{x} \in \widetilde{B}$, so that
$s_d(\gf) \ge l (\gf) \ge 1$. 

\s
Assume finally that $(B,g)$ is a $P$-manifold modelled on $S^2$ 
different from $(S^2, g_\text{can})$, and let $\gt$ be a twist defined
by $g$. 
According to \cite{S0}, $\pi_0 \left( \Symp^c \left(T^*S^2 \right)
\right)$ is generated by the class $[ \tau ]$ of a generalized Dehn
twist $\tau$ defined with respect to $g_{\can}$, and so
$[\gt] = \left[ \tau^k \right]$ for some $k \in \ZZ$.
Clearly, $\gt \neq \id$.
If $k=0$, the estimate $s_1(\gt) = l(\gt) \ge 1$ therefore follows from
Theorem~1, and if $k \neq 0$ from Corollary~3 below.
The proof of Theorem~\ref{t:twist} is complete.
\proofend

\subsection{Proof of Corollary 3}  \label{c:3}

\ni
Let $\gf \in \Sympcc \left( T^*S^d \right)$ be such that
$[\gf] = 
\left[ \tau^m \right] \in \pi_0 \left( \Sympcc \left( T^*S^d \right) \right)$
for some $m \in \ZZ \setminus \{0\}$.
Since $\left[ \tau^2 \right] = [\gt]$, we then have 
$\left[ \gf^2 \right] = \big[ \gt^m \big]$. 
Proceeding as above and assuming again $r=1$ we find $c>0$ such that 
\begin{equation}  \label{est:cn}
\mu_{g^*} \left( \gf^{2n} \left( D_x \right) \right)  \,\ge\, cn 
\end{equation}
for all $n \ge 1$.
We denote by $\left\| D_z \gf \right\|$ 
the operator norm of the differential of $\gf$ at a point $z \in T^*S^d$ 
with respect to $g$, and we abbreviate 
$\left\| D \gf \right\| = \max_{z \in T^*S^d} \left\| D_z \gf \right\|$.
Using the estimate \eqref{est:cn} we find
\begin{eqnarray} \label{est:odd}
\mu_{g^*}  \left( \gf^{2n+1} \left( D_x \right) \right)
&\ge&
\left\| D \gf \right\|^{-1}  
\mu_{g^*}  \left( \gf^{2n+2} \left( D_x \right) \right) \\
&\ge&
\left\| D \gf \right\|^{-1}  c (n+1) . \notag
\end{eqnarray}
The estimates \eqref{est:cn} and \eqref{est:odd} now show that
$l(\gf)\ge 1$, as claimed.
\proofend

\subsection{Differential topology of Dehn twists}  \label{var}

\ni
In this subsection we collect results concerning the differential topology 
of Dehn twists.
We shall in particular see that for odd spheres and their quotients,
Theorem~3 already holds for topological reasons, while for even spheres
and $\CP^n$'s, Theorem~3 is a genuinely symplectic result.

As above, $(B,g)$ is a $d$-dimensional $P$-manifold,
$M = T^*B$ and $M_r = T^*_rB$.
We denote by $\Diffcc (M)$ the group of compactly supported
diffeomorphisms of $M$.
Each $\gf \in \Diffcc (M)$ induces a variation homomorphism
\[
\var_\gf \colon H_*^{cl}(M) \ra\ H_*(M),
\quad  
[ c ] \mapsto [\gf_* c - c] .
\]
Here, the homology $H_*^{cl}(M)$ with closed support as well as $H_*(M)$
are taken with integer coefficients.
Notice that $\gf$ is not isotopic to the identity in $\Diffcc (M)$ if
$\var_\gf \neq 0$.
By Poincar\'e-Lefschetz duality, 
\[
H_*^{cl}(M) \cong H_*^{cl} ( M_r \setminus \pp M_r) 
\cong H_*(M_r, \pp M_r) \cong H^{2d-*}(M_r) \cong H^{2d-*} (B) ,
\] 
and $H_*(M) \cong H_*(B)$, and so
$\var_\gf =0$ except possibly in degree $*=d$.
It is known from classical Picard-Lefschetz theory that
$\var_{\tau} \colon H_d^{cl} \left( T^* S^d \right) \ra H_d \left(
T^*S^d \right)$ does not vanish, see \cite[p.\ 26]{AGV}.
Assume now that $B$ is oriented.
We orient the fibres $T^*_xB$, $x \in B$, such that 
$[B] \cdot \left[ T^*_xB \right] = 1$, 
where $\cdot$ denotes the intersection product in homology determined by
the natural orientation of the cotangent bundle $M$.
Then $H_d^{cl}(M) \cong \ZZ$ is generated by the
fibre class $F = \left[ T^*_xB \right]$, 
and $H_d(M) \cong \ZZ$ is generated by the base class $[B]$, which by
abuse of notation is denoted $B$.

\begin{proposition}  \label{p:var}
Assume that $(B,g)$ is an oriented $P_k$-manifold.
\begin{itemize}
\item[(i)]
If $k$ is even and $d$ is odd, $\var_{\gt^m} (F) = 2 m B$ for $m \in \ZZ$.
\item[(ii)]
If $k$ is odd, $\var_{\gt^m} =0$ for all $m \in \ZZ$.
\end{itemize}
\end{proposition}

\proof
For simplicity we assume again $m \ge 1$.
As in Subsection~\ref{HF} we choose $\gt = \gt_f$, 
fix $x \in B$, choose $y \in W$, and let
$\gt^m \left( T_x^*B \right) \cap T_y^*B  = 
   \big\{ c_0^\pm, \dots, c_{m-1}^\pm \big\}$.
The local intersection number of $\gt^m \left( T_x^*B \right)$ and
$T_y^*B$ at $c_i^\pm$ is $(-1)^{\ind \left( c_i^{\pm} \right)}$.
Recall from Proposition~\ref{p:jacobiper} that $\ind c_0^+ =0$ and 
\[
\ind c_i^- = \ind c_i^+ +k 
\quad
\text{ and }
\quad 
\ind c_{i+1}^+ = \ind c_i^+ + k + d-1 .
\]

\s
\ni
(i)
If $k$ is even and $d$ is odd, we find 
\[
\gt_*^m (F) \cdot F \,=\, 
\sum_{i=0}^{m-1} (-1)^{\ind c_i^+} + (-1)^{\ind c_i^-} \,=\,
\sum_{i=0}^{m-1} 2  \,=\, 2 m,
\]
and so $\var_{\gt^m} (F) \cdot F = \gt_*^m(F) \cdot F - F \cdot F = 2m$,
i.e., $\var_{\gt^m} (F) = 2m B$.

\s
\ni
(ii) 
If $k$ is odd, we find
\[
 \gt_*^m (F) \cdot F \,=\, 
\sum_{i=0}^{m-1} (-1)^{\ind c_i^+} + (-1)^{\ind c_i^-} \,=\, 0,
\]
and so $\var_{\gt^m} (F) \cdot F = 0$, i.e., $\var_{\gt^m} (F) = 0$.
\proofend

Before discussing the variation homomorphism further, let us show how
Proposition~\ref{p:var} (i) leads to a topological proof of Theorem~3
for the known odd-dimensional $P$-manifolds.

\begin{corollary}  \label{c:top}
Assume that $(B,g)$ is a round sphere $S^{2n+1}$ or one of its quotients
$S^{2n+1} / G$ or a Zoll manifold $( S^{2n+1},g )$.
Then the conclusion of Theorem~3 holds true.
Moreover, if $\gf \in \Diffcc \left( T^*B \right)$ is such that 
$[ \gf ] = \left[ \gt^m \right] \in \pi_0 \left( \Diffcc \left( T^*B \right)
\right)$ for some $m \in \ZZ \setminus \{ 0 \}$, then $s_d \left( \gf
\right) \ge 1$. 
\end{corollary}

\proof
According to Proposition~\ref{p:index}, 
$(B,g)$ is a $P_0$-manifold or a $P_{2n}$-manifold, 
and so Proposition~\ref{p:var} (i) shows that
$\gf_*^n (F) = 2mnB +F$ for all $n \ge 1$.
Choose $r < \infty$ so large that $\gf$ is supported in $T_r^*B$, choose
$x \in B$ and set $D_x(r) = T_x^*B \cap T_r^*B$. Then
\[
\mu_{g^*} \left( \gf^n \left( D_x(r) \right) \right) \,\ge\, 
\left( 2m \mu_g (B) \right) n
\]
for all $n \ge 1$, and so the corollary follows.
\proofend

We now turn to the known even-dimensional $P$-manifolds.
Propositions~\ref{p:index} and \ref{p:var} (ii) show that 
$\var_{\gt^m} = 0$ for 
$S^{2n}$, $\CP^n$, $\HP^n$, $\Ca$ and even-dimensional Zoll manifolds
and all $m \in \ZZ$.
For the non-orientable spaces $\RP^{2n}$ and $\CP^{2n-1} / \ZZ_2$ the
vanishing of $\var_{\gt^m}$ follows from 
$H_{2n} \left( \RP^{2n} \right) = 0$ and 
$H_{4n-2} \left( \CP^{2n-1} / \ZZ_2 \right) = 0$.
The variation homomorphism can be defined for homology with coefficients
in any Abelian group $G$, and one checks that
$\var_{\gt^m}$ vanishes over any finitely generated $G$ for all the
above even-dimensional $P$-manifolds and every $m \in \ZZ$.
Note that if $\var_{\gt^m} \neq 0$ for some $m \neq 0$
then $\gt$ is not isotopic to the
identity in $\Diffcc (M)$. Since we are not aware of another
obstruction we ask

\begin{question}  \label{q:isotopic}
{\rm
Assume that $(B,g)$ is one of $\RP^{2n}$, $\HP^n$, $\Ca$, $\CP^{2n-1} /
\ZZ_2$.
Is $\gt$ isotopic to the identity in $\Diffcc \left( T^*B \right)$?
}
\end{question}

We did not ask Question \ref{q:isotopic} for 
even-dimensional Zoll manifolds or $\CP^n$ in view of

\begin{proposition}  \label{p:iso}
Assume that $(B,g)$ is an even-dimensional Zoll manifold or $\CP^n$.
Then $\gt$ is isotopic to the identity in $\Diffcc \left( T^*B \right)$. 
\end{proposition} 

\proof
This has been proved in \cite{S2} for $B = \CP^n$, $n \ge 1$, 
by extending the construction for $S^2$ given in \cite{S1}.
This construction carries over literally to $S^6$ since $S^6$ carries an
almost complex structure induced by the vector product on $\RR^7$
related to the Cayley numbers, see \cite[Example~4.4]{MS}.
For arbitrary even-dimensional spheres the result is proved by 
N.\ Krylov, \cite{K}.
Finally, let $\gt = \gt_f$ be a twist defined by a Zoll metric $g$ on $S^{2n}$.
Then there is a smooth family $g_t$, $t \in [0,1]$, of $P$-metrics
with $g_0 = g_\text{can}$ and $g_1 = g$.
The family $\gt_t$ of twists defined by $f$ and $g_t$ is then a
smooth family in $\Sympcc \left( T^*S^{2n} \right)$. 
In particular, $\gt = \gt_1$ is isotopic to $\gt_0$ and hence to the
identity in $\Diff^c \left( T^*S^{2n} \right)$.
\proofend

\begin{remark}  \label{r:ac}
{\rm
One can show that $\gt$ is isotopic to the identity in 
$\Diffcc \left( T^*B \right)$
for any {\it almost complex} $P$-manifold $B$. 
The only almost complex manifolds among the known $P$-manifolds are,
however, $\CP^n$, $n \ge 1$, and $S^6$.
}
\end{remark}

\section{Remarks on smoothness}  \label{smooth}

\ni
Entropy type estimates often depend on the differentiability of the maps
considered. E.g., in the case of finite smoothness the entropy
conjecture \eqref{est:ent} has been proved only for special classes of
manifolds and maps, and there exist homeomorphisms $\gf$ of compact
manifolds such that $h_\text{top}(\gf) < \log \rho (\gf_*)$, see \cite{K1}.   
The results established in this paper hold under essentially minimal
differentiability assumptions necessary to formulate them.
This is so because the uniform lower bounds found for $s_1(\gf)$,
$s_d(\gf)$ and $l(\gf)$ and $C^{\infty}$-smooth symplectomorphisms are
of ``symplecto-topological'' nature.
Given a symplectic manifold $(M,\go)$, let $\Symp^{c,1} (M,\go)$ be
the group of $C^1$-smooth symplectomorphisms whose support is
compact and contained in $M \setminus \pp M$,
and let $\Symp^{c,1}_0 (M,\go)$ be its identity component.
\begin{proposition}
Theorem~\ref{t:4exact} holds true for $\gf \in \Symp^{c,1}_0 (M,d\gl)$.
\end{proposition}

\proof
Proposition~\ref{p:4axa} is proved in \cite{FS} for $C^2$-smooth
Hamiltonians. The Flux is defined on $\Symp^{c,1}_0 (M,\go)$, 
and the exact sequence \eqref{s:exact} exists in the 
$C^1$-category. The remaining arguments in the proof of
Theorem~\ref{t:4exact} are of topological nature and thus go
through for $C^1$-smooth symplectomorphisms.
\proofend

Corollary~\ref{c:torsion}, Theorem~\ref{t:distortion} and
Corollary~\ref{c:finite} also continue to hold in the $C^1$-category.
\begin{proposition}
Theorem~2 holds true for $C^2$-smooth classical Hamiltonian functions.
\end{proposition}

\proof
The results of \cite{Be} hold for $C^2$-smooth Hamiltonians, and
the arguments in Section~\ref{t2} go through.
\proofend

We endow $\Symp^{c,1} (M, \go)$ with the
$C^1$-topology. According to a result of Zehnder, \cite{Z}, 
$\Symp^c (M, \go)$ is dense in $\Symp^{c,1} (M, \go)$, and by a
result of Weinstein, \cite[Theorem~10.1]{MS}, both groups are
locally path connected. It follows that the inclusion 
$\Symp^c (M, \go) \ra \Symp^{c,1} (M, \go)$ induces an
isomorphism of mapping class groups, 
$\pi_0 \left( \Symp^c (M, \go) \right) = 
\pi_0 \left( \Symp^{c,1} (M, \go) \right)$.

\begin{proposition}
Theorem~\ref{t:twist} and Corollary~3 hold true for $C^1$-smooth
symplectomorphisms.
\end{proposition}

\proof
Let $(B,g)$ be as in Theorem~\ref{t:twist}, and let 
$\gf \in \Symp^c \left( T^*B \right)$ be a $C^1$-smooth
symplectomorphism such that $[ \gf ] = \left[ \gt^m \right] \in \pi_0
\left( \Symp^c \left(T^*B \right) \right)$ for some 
$m \in \ZZ \setminus \{ 0 \}$.
We can assume that $\gf$ is supported in $T_1^*B$.
Choose a sequence $\gf_i \in \Symp^c \left(T^*_1B \right)$
of $C^\infty$-smooth symplectomorphisms such that $\gf_i \ra \gf$ in the
$C^1$-topology.
For $i$ large enough, 
$[\gf_i] = [\gf] \in \pi_0 \left( \Symp^c \left(T^*B \right) \right)$.
Using the estimate \eqref{est:fU} we thus conclude 
\[
\mu_{g^*} \big( \gf^n (D_x) \big) \,=\,
\lim_{i \ra \infty} \mu_{g^*} \big( \gf_i^n (D_x) \big) \,\ge\, 
2mn \mu_g(U) 
\]
for all $n \ge 1$. Therefore, $s_d(\gf) \ge l (\gf) \ge
1$. Corollary~3 now follows also for $C^1$-smooth symplectomorphisms.
\proofend

\section{Comparison of volume growth and growth of the differential}
\label{comp}

\ni
For any compactly supported diffeomorphism $\gf$ of a smooth manifold $M$
we denote by $\left\| D_x\gf \right\|$ 
the operator norm of the differential of $\gf$ at the point $x$ 
with respect to some Riemannian metric on $M$.
Following \cite{DG} we define the {\it growth sequence}\, of $\gf$ by 
\[
\Gamma_n(\gf) \,=\, \max_{x\in M} \left\| D_x \gf^n \right\| .
\]
In \cite{P1}, Polterovich proved that if $(M, \go)$
is a closed symplectic manifold with $\pi_2(M) = 0$, then 
for a large class of symplectomorphisms $\gf \in \Symp_0(M, \go)$
there exists a uniform lower bound for the growth type of the sequence 
$\Gamma_n(\gf)$.
Complementary results for symplectic and smooth diffeomorphisms were
found in \cite{P2} and \cite{PSo}.
Here, we derive from the growth sequence the {\it slow differential growth}
\begin{equation*} 
\gg (\gf ) \,=\, \liminf_{n \ra \infty} \frac{\log \Gamma_n(\gf)}{\log n}.
\end{equation*}
It does not depend on the choice of Riemannian metric.
While in the definition of the slow volume growth $s_i(\gf)$ we looked
at the weighted asymptotics of the most distorted smooth $i$-dimensional
family  
of orbits, in the definition of $\gg (\gf)$ we look at each time $n$ at
the place of the strongest distortion and pass to a weighted limit.
The invariants $s_i(\gf)$ are thus of rather dynamical nature, while
$\gg(\gf)$ is of rather geometric nature.
Clearly,
\begin{equation} \label{est:sig} 
s_i (\gf) \,\le\, i \,\gg(\gf),  \quad \, i = 1, \dots, \dim M .
\end{equation} 
We therefore read off from our main results
\begin{corollary}
Under the assumptions of Theorem \ref{t:4exact}
it holds that $\gg (\gf) \ge 1$.
\end{corollary}
\begin{corollary}
Under the assumptions of Theorem~2
it holds that $\gg (\gf) \ge 1/2$ respectively $\gg (\gf) \ge 1$.
\end{corollary}
\begin{corollary}  \label{c:1k}
Under the assumptions of Theorem~\ref{t:twist} or Corollary~3, 
it holds that $\gg (\gf) \ge 1/d$.
\end{corollary}

Our proof of Proposition~\ref{prop:n} (ii) shows that $\gg (\gt) =1$ for
any twist $\gt$ of the cotangent $T^*B$ over a $P$-manifold $(B,g)$.
It follows that $\gg (\tau) = 1$ for any (generalized) Dehn twist $\tau$
of $T^*S^d$.

\begin{question}  \label{q:sharp}
{\rm
Can the lower bound $1/d$ in Corollary \ref{c:1k} be replaced by $1$? 
}
\end{question}

The estimates \eqref{est:sig} are in general not sharp:
%
%
%
Choose a monotone decreasing smooth function $f \colon \RR \to
[1,2]$ with 
$f(r) = 2$ if $r \le 1$ and $f(r) =1$ if $r \ge 2$, and define
$\phi \in \Diff_0^c \left( \RR^n \right)$ by
$\phi (x) = f \left( |x| \right) x$. 
Using $\phi$ as a plug, we see that every smooth manifold $M$
carries $\gf \in \Diff_0^c (M)$ with $s(\gf)=0$ and $\gg (\gf) =
\infty$.
It would be interesting to find such diffeomorphisms in the
symplectic category.

\enddocument